\documentclass[12pt]{article}%
\usepackage{amsmath}
\usepackage{amsfonts}
\usepackage{amssymb}
\usepackage{graphicx}%
\providecommand{\U}[1]{\protect\rule{.1in}{.1in}}
\newtheorem{theorem}{Theorem}[section]

\newtheorem{claim}[theorem]{Claim}

\newtheorem{definition}{Definition}[section]

\newtheorem{lemma}[theorem]{Lemma}

\newtheorem{observation}[theorem]{Observation}

\newtheorem{proposition}[theorem]{Proposition}
\newtheorem{question}{Question}

\providecommand{\boksie}{\ensuremath{\mathbin{\raisebox{0.3mm}{$\scriptstyle\square$}}}}
\begin{document}

\title{\textbf{Comparing Upper Broadcast Domination and Boundary Independence
Broadcast Numbers of Graphs}}
\author{C.M. Mynhardt\thanks{Supported by the Natural Sciences and Engineering
Research Council of Canada.}\\Department of Mathematics and Statistics\\University of Victoria, Victoria, BC, \textsc{Canada}\\{\small kieka@uvic.ca}
\and L. Neilson\\Department of Adult Basic Education\\Vancouver Island University, Nanaimo, BC, \textsc{Canada}\\{\small linda.neilson@viu.ca}}
\maketitle

\begin{abstract}
A broadcast on a nontrivial connected graph $G=(V,E)$ is a function
$f:V\rightarrow\{0,1,\dots,\operatorname{diam}(G)\}$ such that $f(v)\leq e(v)$
(the eccentricity of $v$) for all $v\in V$. The weight of $f$ is $\sigma(f)=%
{\textstyle\sum_{v\in V}}
f(v)$. A vertex $u$ hears $f$ from $v$ if $f(v)>0$ and $d(u,v)\leq f(v)$. A
broadcast $f$ is dominating if every vertex of $G$ hears $f$. The upper
broadcast domination number of $G$ is $\Gamma_{b}(G)=\max\left\{
\sigma(f):f\text{ is a minimal dominating broadcast of }G\right\}  .$

A broadcast $f$ is boundary independent if, for any vertex $w$ that hears $f$
from vertices $v_{1},...,v_{k},\ k\geq2$, the distance $d(w,v_{i})=f(v_{i})$
for each $i$. The maximum weight of a boundary independent broadcast is the
boundary independence broadcast number $\alpha_{\operatorname{bn}}(G)$.

We compare $\alpha_{\operatorname{bn}}$ to $\Gamma_{b}$, showing that neither
is an upper bound for the other. We show that the differences $\Gamma
_{b}-\alpha_{\operatorname{bn}}$ and $\alpha_{\operatorname{bn}}-\Gamma_{b}$
are unbounded, the ratio $\alpha_{\operatorname{bn}}/\Gamma_{b}$ is bounded
for all graphs, and $\Gamma_{b}/\alpha_{\operatorname{bn}}$ is bounded for
bipartite graphs but unbounded in general.

\end{abstract}

\noindent\textbf{Keywords:\hspace{0.1in}}broadcast domination; broadcast
independence; hearing independent broadcast; boundary independent broadcast;
boundary independence broadcast number

\noindent\textbf{AMS Subject Classification Number 2010:\hspace{0.1in}}05C69

\label{Submitted to Trans. Comb., March 23, 2021}

\section{Introduction}

The study of broadcast domination and broadcast independence was initiated by
Erwin in his doctoral dissertation \cite{Ethesis}. To generalize the concept
of an independent set $X$ in a graph $G$ to an independent broadcast, Erwin
focussed on the property that no vertex in $X$ belongs to the neighbourhood of
another vertex in $X$. Focussing instead on the property that no edge of $G$
is incident with (or covered by) more than one vertex in $X$, Neilson
\cite{LindaD} and Mynhardt and Neilson \cite{MN} defined boundary independent
broadcasts as an alternative to Erwin's independent broadcasts. We explain
below why this definition results in a parameter, called the boundary
independence number $\alpha_{\operatorname{bn}}(G)$, that is, in some sense,
\textquotedblleft better behaved\textquotedblright\ than Erwin's broadcast
independence number, which we denote here by $\alpha_{h}(G)$.

We compare the upper broadcast domination number $\Gamma_{b}$, also defined by
Erwin \cite{Ethesis}, to $\alpha_{\operatorname{bn}}$, showing that neither is
an upper bound for the other. We denote this incomparability by $\alpha
_{\operatorname{bn}}\diamond\Gamma_{b}$. We show that the differences
$\Gamma_{b}-\alpha_{\operatorname{bn}}$ and $\alpha_{\operatorname{bn}}%
-\Gamma_{b}$ are unbounded, the ratio $\alpha_{\operatorname{bn}}/\Gamma_{b}$
is bounded for all graphs, and $\Gamma_{b}/\alpha_{\operatorname{bn}}$ is
bounded for bipartite graphs but unbounded in general.

\section{Definitions and background}

\label{SecDefs}For undefined concepts we refer the reader to \cite{CLZ}. A
\emph{broadcast} on a connected graph $G=(V,E)$ is a function $f:V\rightarrow
\{0,1,\dots,\operatorname{diam}(G)\}$ such that $f(v)\leq e(v)$ (the
eccentricity of $v$) for all $v\in V$ if $|V|\geq2$, and $f(v)=1$ if
$V=\{v\}$. When $G$ is disconnected, we define a broadcast on $G$ as the union
of broadcasts on its components. Define $V_{f}^{+}=\{v\in V:f(v)>0\}$ and
partition $V_{f}^{+}$ into the two sets $V_{f}^{1}=\{v\in V:f(v)=1\}$ and
$V_{f}^{++}=V_{f}^{+}-V_{f}^{1}$. A vertex in $V_{f}^{+}$ is called a
\emph{broadcasting vertex}. A vertex $u$ \emph{hears} $f$ from $v\in V_{f}%
^{+}$, and $v$ $f$-\emph{dominates} $u$, if the distance $d(u,v)\leq f(v)$.
Denote the set of all vertices that do not hear $f$ by $U_{f}$. A broadcast
$f$ is \emph{dominating} if $U_{f}=\varnothing$. For any subset $U$ of $V$, we
define $f(U)=\sum_{u\in U}f(u)$. The \emph{weight} of $f$ is $\sigma(f)=f(V)$,
and the \emph{broadcast number} of $G$ is
\[
\gamma_{b}(G)=\min\left\{  \sigma(f):f\text{ is a dominating broadcast of
}G\right\}  .
\]
If $f$ and $g$ are broadcasts on $G$ such that $g(v)\leq f(v)$ for each $v\in
V$, we write $g\leq f$. If in addition $g(v)<f(v)$ for at least one $v\in V$,
we write $g<f$. A dominating broadcast $f$ on $G$ is a \emph{minimal
dominating broadcast} if no broadcast $g<f$ is dominating. The \emph{upper
broadcast domination number }of $G$ is
\[
\Gamma_{b}(G)=\max\left\{  \sigma(f):f\text{ is a minimal dominating broadcast
of }G\right\}  ,
\]
and a dominating broadcast $f$ of $G$ such that $\sigma(f)=\Gamma_{b}(G)$ is
called a $\Gamma_{b}(G)$-\emph{broadcast} (abbreviated to $\Gamma_{b}%
$-\emph{broadcast }if the graph $G$ is obvious). Introduced by Erwin
\cite{Ethesis}, the upper broadcast number was also studied by, for example,
Ahmadi, Fricke, Schroeder, Hedetniemi and Laskar \cite{Ahmadi}, Bouchemakh and
Fergani \cite{BF}, Bouchouika, Bouchemakh and Sopena \cite{BBS}, Dunbar,
Erwin, Haynes, Hedetniemi and Hedetniemi \cite{Dunbar}, and Mynhardt and Roux
\cite{MR}.

We denote the independence number of $G$ by $\alpha(G)$; an independent set of
$G$ of cardinality $\alpha(G)$ is called an $\alpha(G)$\emph{-set}, often
abbreviated to $\alpha$-\emph{set}. To generalize the concept of independent
sets, Erwin \cite{Ethesis} defined a broadcast $f$ to be \emph{independent},
or, for our purposes, \emph{hearing independent}, abbreviated to
\emph{h-independent}, if no vertex $u\in V_{f}^{+}$ hears $f$ from any other
vertex $v\in V_{f}^{+}$; that is, broadcasting vertices only hear themselves.
The maximum weight of an h-independent broadcast is the \emph{h-independent
broadcast number}, denoted by $\alpha_{h}(G)$; such a broadcast is called an
$\alpha_{h}(G)$\emph{-broadcast }(or $\alpha_{h}$\emph{-broadcast} for short).
This version of broadcast independence was also considered by, among others,
Ahmane, Bouchemakh and Sopena \cite{ABS, ABS2}, Bessy and Rautenbach \cite{BR,
BR2}, Bouchemakh and Zemir \cite{Bouch}, Bouchouika et al.~\cite{BBS} and
Dunbar et al.~\cite{Dunbar}. For a survey of broadcasts in graphs, see the
chapter by Henning, MacGillivray and Yang \cite{HMY}.

Before continuing, we define a class of trees often used as examples. For
$k\geq3$ and $n_{i}\geq1$ for\ $i\in\{1,...,k\}$, the \emph{(generalized)
spider} $\operatorname{Sp}(n_{1},...,n_{k})$ is the tree which has exactly one
vertex $b$, called the \emph{head}, having $\deg(b)=k$, and for which the $k$
components of $\operatorname{Sp}(n_{1},...,n_{k})-b$ are paths of lengths
$n_{1}-1,...,n_{k}-1$, respectively. The \emph{legs }$L_{1},...,L_{k}$ of the
spider are the paths from $b$ to the leaves. If $n_{i}=r$ for each $i$, we
write $\operatorname{Sp}(r^{k})$ for $\operatorname{Sp}(n_{1},...,n_{k})$. An
\emph{endpath} in a tree is a path ending in a leaf and having all internal
vertices (if any) of degree $2$; the legs of a spider are examples of endpaths.

\subsection{Neighbourhoods and boundaries}

\label{SecNBHs}For a broadcast $f$ on a graph $G$ and $v\in V_{f}^{+}$, we
define the%
\begin{equation}
\left.
\begin{tabular}
[c]{l}%
$f$-\emph{neighbourhood}\\
$f$-\emph{boundary}\\
$f$-\emph{private }\\%
\begin{tabular}
[c]{l}%
\
\end{tabular}
\emph{neighbourhood}\\
\\
$f$-\emph{private boundary}\\%
\begin{tabular}
[c]{l}%
\
\end{tabular}
\\%
\begin{tabular}
[c]{l}%
\
\end{tabular}
\end{tabular}
\ \right\}  ~\text{of }v\text{ by~}\left\{  \text{%
\begin{tabular}
[c]{rll}%
$N_{f}(v)$ & $=$ & $\{u\in V:d(u,v)\leq f(v)\}$\\
$B_{f}(v)$ & $=$ & $\{u\in V:d(u,v)=f(v)\}$\\
&  & \\
$\operatorname{PN}_{f}(v)$ & $=$ & $\{u\in N_{f}(v):u\notin N_{f}(w)$\\
&  & for all$\ w\in V_{f}^{+}-\{v\}\}$\\
$\operatorname{PB}_{f}(v)$ & $=$ & $\{u\in N_{f}(v):u$ is not\\
&  & dominated by$\ (f-\{(v,f(v)\})$\\
&  & $\cup\{(v,f(v)-1)\}$.
\end{tabular}
}\right.  \label{eq_sets_defs}%
\end{equation}

If $v\in V_{f}^{1}$ and $v$ does not hear $f$ from any vertex $u\in V_{f}%
^{+}-\{v\}$, then $v\in\operatorname{PB}_{f}(v)$, and if $v\in V_{f}^{++}$,
then $\operatorname{PB}_{f}(v)=B_{f}(v)\cap\operatorname{PN}_{f}(v)$. Since
$f(v)\leq e(v)$ for each vertex $v$ of a non-trivial graph $G$, $B_{f}%
(v)\neq\varnothing$. Therefore, if there exists a vertex $v$ such that
$\operatorname{PB}_{f}(v)=\varnothing$, then each vertex in $B_{f}(v)$ also
hears $f$ from another broadcasting vertex and we deduce that $|V_{f}^{+}%
|\geq2$. If $f$ is a broadcast such that every vertex $x$ that hears more than
one broadcasting vertex also satisfies $d(x,u)\geq f(u)$ for all $u\in
V_{f}^{+}$, we say that the \emph{broadcast overlaps only in boundaries}.
Equivalently, $f$ overlaps only in boundaries if $N_{f}(u)\cap N_{f}%
(v)\subseteq B_{f}(u)\cap B_{f}(v)$ for all distinct $u,v\in V_{f}^{+}$. If
$uv\in E(G)$ and $u,v\in N_{f}(x)$ for some $x\in V_{f}^{+}$ such that at
least one of $u$ and $v$ does not belong to $B_{f}(x)$, we say that the edge
$uv$ is \emph{covered} in $f$, or $f$-\emph{covered}, by $x$. If $uv$ is not
covered by any $x\in V_{f}^{+}$, we say that $uv$ is \emph{uncovered by~}$f$
or $f$-\emph{uncovered}. We denote the set of $f$-uncovered edges by
$U_{f}^{E}$.

Erwin \cite{Ethesis} determined a necessary and sufficient condition for a
dominating broadcast to be minimal dominating. We restate it here in terms of
private boundaries.

\begin{proposition}
\label{PropMinimal}\emph{\cite{Ethesis}}\hspace{0.1in}A dominating broadcast
$f$ is a minimal dominating broadcast if and only if $\operatorname{PB}%
_{f}(v)\neq\varnothing$ for each $v\in V_{f}^{+}$.
\end{proposition}

Ahmadi et al.~\cite{Ahmadi} define \label{Def_ir}a broadcast $f$ to be
\emph{irredundant} if $\operatorname{PB}_{f}(v)\neq\varnothing$ for each $v\in
V_{f}^{+}$. Therefore, a dominating broadcast is minimal dominating if and
only if it is irredundant.

\subsection{Boundary independent broadcasts}

\label{bn-bcs} The characteristic function of an independent set of a graph
$G$ is a broadcast on $G$ that overlaps only in boundaries. This feature was
used in \cite{MN, LindaD} to define three types of \emph{boundary independent
broadcasts}; since no edge is covered by more than one broadcasting vertex,
these types of broadcasts are, in some sense, more independent than
h-independent broadcasts.

\begin{definition}
\label{bn-i}\emph{\cite{MN, LindaD}\hspace{0.1in}A broadcast is}
bn-independent \emph{if it overlaps only in boundaries. The maximum weight of
a bn-independent broadcast on }$G$\emph{ is }$\alpha_{\operatorname{bn}}%
(G)$\emph{; such a broadcast is called an }$\alpha_{\operatorname{bn}}%
(G)$-broadcast\emph{.}
\end{definition}%

\begin{figure}[ptb]%
\centering
\includegraphics[
height=2.4898in,
width=4.1079in
]%
{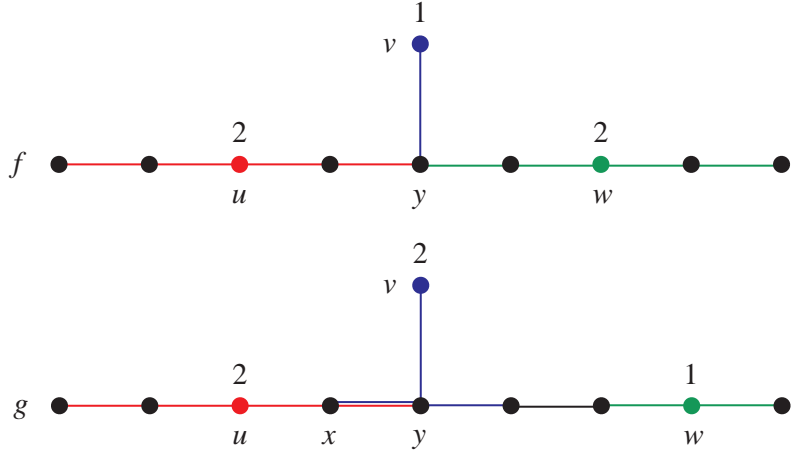}%
\caption{The broadcast $f$ (top) is bn-independent, but $g$ is not}%
\label{Fig_def2_1}%
\end{figure}
For example, the broadcast $f$ in Figure \ref{Fig_def2_1} is bn-independent,
because $y$ is the only vertex that hears $f$ from more than one vertex in
$V_{f}^{+}$, and $y\in B_{f}(u)\cap B_{f}(v)\cap B_{f}(w)$. On the other hand,
the broadcast $g$ is not boundary independent because $x$ hears $g$ from $u$
and $v$, but $x\notin B_{g}(u)$ (and $y\notin B_{g}(v)$); the edge $xy$ is
covered by $u$ and by $v$.

Note that in a bn-independent broadcast, no broadcasting vertex $u$ hears
another broadcasting vertex $v$ (otherwise an edge incident with $u$ or $v$
would be covered by both $u$ and $v$), hence any bn-independent broadcast is
also hearing independent.

\begin{definition}
\label{bnr-i}\emph{\cite{MN, LindaD}\hspace{0.1in}A broadcast is}
bnr-independent\emph{ if it is bn-independent and irredundant. The maximum
weight of a bnr-independent broadcast is} $\alpha_{\operatorname{bnr}}%
(G)$\emph{; such a broadcast is called an }$\alpha_{\operatorname{bnr}}(G)$-broadcast
\end{definition}

\begin{definition}
\label{bnd-i}\emph{\cite{MN, LindaD}\hspace{0.1in}A broadcast is}
bnd-independent\emph{ if it is minimal dominating and bn-independent. The
maximum weight of a bnd-independent broadcast is} $\alpha_{\operatorname{bnd}%
}(G)$\emph{; such a broadcast is called an }$\alpha_{\operatorname{bnd}}%
(G)$-broadcast\emph{.}
\end{definition}

When the graph $G$ is clear, an $\alpha_{\operatorname{bn}}(G)$-broadcast is
also called an $\alpha_{\operatorname{bn}}$-\emph{broadcast}; the same comment
holds for $\alpha_{\operatorname{bnr}}(G)$ and $\alpha_{\operatorname{bnd}%
}(G)$. Since the characteristic function of an independent set is a bnd-,
bnr-, bn- and h-independent broadcast, it follows from Definitions \ref{bn-i}
-- \ref{bnd-i} that%
\begin{equation}
\alpha(G)\leq\alpha_{\operatorname{bnd}}(G)\leq\alpha_{\operatorname{bnr}%
}(G)\leq\alpha_{\operatorname{bn}}(G)\leq\alpha_{h}(G)\label{sequence2}%
\end{equation}
for any graph $G$. As shown in \cite{MN, LindaD}, if $G$ is a $2$-connected
bipartite graph, then
\begin{equation}
\alpha(G)=\alpha_{\operatorname{bnd}}(G)=\alpha_{\operatorname{bnr}}%
(G)=\alpha_{\operatorname{bn}}(G).\label{sequence5}%
\end{equation}
Denote the $m\times n$ grid by $G_{m,n}$. Bouchemakh and Zemir \cite{Bouch}
showed that $\alpha_{h}(G_{5,5})=15>\left\lceil \frac{25}{2}\right\rceil
=\alpha(G_{5,5})$, hence the equality string (\ref{sequence5}) does not extend
to $\alpha_{h}$. If $v$ is a peripheral vertex of a connected graph $G\neq
K_{1}$, then broadcasting with a strength of $\operatorname{diam}(G)$ from $v$
and $0$ from any other vertex results in a minimal dominating bn-independent
broadcast, hence we also have
\[
\alpha_{\operatorname{bn}}(G)\geq\operatorname{diam}(G)\text{ and }\Gamma
_{b}(G)\geq\operatorname{diam}(G)
\]
for all nontrivial connected graphs $G$.

Suppose $f$ is a bn-, bnr- or bnd-independent broadcast on $G$ and an edge
$uv$ of $G$ is covered by vertices $x,y\in V_{f}^{+}$. By the definition of
covered, $\{u,v\}\nsubseteq B_{f}(x)$ and $\{u,v\}\subseteq N_{f}(x)\cap
N_{f}(y)$. This violates the bn-independence of $f$. Hence we have the
following property, a primary motivation for the definition of bn-independent broadcasts.

\begin{observation}
\label{edge-covered}If $f$ is a bn-, bnr- or bnd-independent broadcast on a
graph $G$, then each edge of $G$ is covered by at most one vertex in
$V_{f}^{+}$.
\end{observation}

Bouchemakh and Fergani \cite{BF} showed that if $G$ is a graph of order $n$
and minimum degree $\delta(G)$, then $\Gamma_{b}(G)\leq n-\delta(G)$; the
bound is sharp, for example for paths, stars and complete graphs. Mynhardt and
Neilson \cite{MN} showed that $\alpha_{\operatorname{bn}}(G)\leq n-1$ for all
graphs $G$, and that equality holds for a connected graph $G$ if and only if
$G$ is a path or a generalized spider. Hence $\Gamma_{b}(G)$ and
$\alpha_{\operatorname{bn}}(G)$ never exceed the order of $G$. On the other
hand, it follows from a result in \cite{Dunbar} that $\alpha_{h}%
(\operatorname{Sp}(r^{k}))=k(2r-1)$. Since $\operatorname{Sp}(r^{k})$ has
order $kr+1$, we see that there exist graphs whose h-independent broadcast
number is almost double their order. One may thus regard $\alpha
_{\operatorname{bn}}$ as being \textquotedblleft better
behaved\textquotedblright\ than $\alpha_{h}$.

A bn-independent broadcast $f$ on $G$ is \emph{maximal bn-independent} if no
broadcast $g$ on $G$ such that $g>f$ is bn-independent. The minimum weight of
a maximal bn-independent broadcast on $G$ is the \emph{lower broadcast
independence number}, denoted by $i_{\operatorname{bn}}(G)$. This parameter
was introduced by Neilson \cite{LindaD} and further investigated by
Marchessault and Mynhardt \cite{MM}.

The rest of the paper is organized as follows. We present previous results in
Section \ref{Sec_known}. In Section \ref{Sec_Non_Comp} we show that
$\alpha_{\operatorname{bnr}}\diamond\Gamma_{b}$ and $\alpha_{\operatorname{bn}%
}\diamond\Gamma_{b}$, that is, neither pair of parameters are comparable. In
Section \ref{SecDiff} we show that the differences $\Gamma_{b}-\alpha
_{\operatorname{bn}}$ and $\Gamma_{b}-\alpha_{\operatorname{bnr}}$ can be
arbitrary for general graphs, while $\Gamma_{b}-\alpha_{\operatorname{bnr}}$
can also be arbitrary for trees. In the other direction we show that
$\alpha_{\operatorname{bnr}}-\Gamma_{b}$, hence also $\alpha
_{\operatorname{bn}}-\Gamma_{b}$, can be arbitrary for trees as well as for
cyclic graphs. We consider the ratios of these parameters in Section
\ref{SecRatios}, showing that $\alpha_{\operatorname{bnr}}(G)/\Gamma
_{b}(G)\leq\alpha_{\operatorname{bn}}(G)/\Gamma_{b}(G)<2$ for all graphs,
whereas $\Gamma_{b}(G)/\alpha_{\operatorname{bnr}}(G)$ and $\Gamma
_{b}(G)/\alpha_{\operatorname{bn}}(G)$ are unbounded for general graphs, but
$\Gamma_{b}(G)/\alpha_{\operatorname{bn}}(G)\leq\Gamma_{b}(G)/\alpha
_{\operatorname{bnr}}(G)<2$ for connected bipartite graphs. We conclude with
questions for future consideration in Section \ref{Sec_Qs}.

\subsection{Known results}

\label{Sec_known}In this subsection we present known results that will be used
later on. It is often useful to know when a bn-independent broadcast is
maximal bn-independent.

\begin{proposition}
\label{prop-max-bn}

\begin{enumerate}
\item[$(i)$] \emph{\cite{MN}}\hspace{0.1in}A bn-independent broadcast $f$ on a
graph $G$ is maximal bn-indepen-dent if and only if (a) it is dominating, and
(b) either $V_{f}^{+}=\{v\}$ or $B_{f}(v)-\operatorname{PB}_{f}(v)\neq
\varnothing$ for each $v\in V_{f}^{+}$.

\item[$(ii)$] \emph{\cite{MM}\hspace{0.1in}}Let $f$ be a bn-independent
broadcast on a connected graph $G$ such that $|V_{f}^{+}|\geq2$. Then $f$ is
maximal bn-independent if and only if each component of $G-U_{f}^{E}$ contains
at least two broadcasting vertices.
\end{enumerate}
\end{proposition}

When $f$ is a minimal dominating or a bnr-broadcast, $\operatorname{PB}%
_{f}(v)\neq\varnothing$ for each $v\in V_{f}^{+}$, but when $f$ is an h- or
bn-independent broadcast, it is possible that $\operatorname{PB}%
_{f}(v)=\varnothing$. The next observations follow from (\ref{eq_sets_defs})
in Section \ref{SecNBHs} and Definitions \ref{bn-i} -- \ref{bnd-i} in Section
\ref{bn-bcs}; we state them here for referencing.

\begin{observation}
\label{Ob_PB}\emph{\cite{MN, LindaD}\hspace{0.1in}}$(i)\hspace{0.1in}$If $f$
is a bn-, bnr- or h-independent broadcast and $v\in V_{f}^{1}$, then
$v\in\operatorname{PB}_{f}(v)$.

\begin{enumerate}
\item[$(ii)$] If $f$ is an h- or a bn-independent broadcast such that
$\operatorname{PB}_{f}(v)=\varnothing$ for some $v\in V_{f}^{+}$, then $v\in
V_{f}^{++}$ and $|V_{f}^{+}|\geq2$.
\end{enumerate}
\end{observation}

Let $f$ be an $\alpha_{\operatorname{bn}}$-broadcast on a graph $G$ and let
$T$ be a spanning tree of $G$. Removing the edges in $E(G)-E(T)$ does not
affect bn-independence, hence $f$ is also a bn-independent broadcast on $T$.
This leads to the following observation.

\begin{observation}
\label{spanning_tree}If $T$ is a spanning tree of a graph $G$, then
$\alpha_{\operatorname{bn}}(T)\geq\alpha_{\operatorname{bn}}(G)$.
\end{observation}

The above inequality can be strict, for example for $K_{n},\ n\geq3$, because
$\alpha_{\operatorname{bn}}(K_{n})=1$ and $\alpha_{\operatorname{bn}}%
(T)\geq\operatorname{diam}(T)\geq2$ for all trees of order at least $3$. 

When determining boundary independence numbers of trees, the following results
from \cite{MN2, LindaD}, which we summarize in a single theorem, can be quite useful.

\begin{theorem}
\label{ThmTrees}\hspace{0.1in}Consider a tree $T$.

\begin{enumerate}
\item[$(i)$] \emph{\cite{MN2}}, \emph{\cite[Lemma 2.3.12]{LindaD}%
\hspace{0.1in}}If $f$ is an $\alpha_{\operatorname{bn}}(T)$-broadcast, then no
leaf hears $f$ from a non-leaf.

\item[$(ii)$] \emph{\cite[Theorem 2.3.13]{LindaD}\hspace{0.1in}}There exists
an $\alpha_{\operatorname{bnr}}(T)$-broadcast $f$ such that no leaf hears $f$
from a non-leaf.

\item[$(iii)$] \emph{\cite{MN2}}, \emph{\cite[Theorem 2.3.14]{LindaD}%
\hspace{0.1in}}There exists an $\alpha_{\operatorname{bn}}(T)$-broadcast $f$
such that for all $v\in V^{+}$, $f(v)=1$ or $\deg(v)=1$.

\item[$(iv)$] \emph{\cite[Theorem 2.3.14]{LindaD}\hspace{0.1in}}There exists
an $\alpha_{\operatorname{bnr}}(T)$-broadcast $f$ such that for all $v\in
V^{+}$, $f(v)=1$ or $\deg(v)=1$.

\item[$(v)$] \emph{\cite{MN2}}, \emph{\cite[Lemma 2.3.15]{LindaD}%
\hspace{0.1in}}If $f$ is an $\alpha_{\operatorname{bn}}(T)$-broadcast such
that $|V_{f}^{1}|$ is maximum, then $\operatorname{PB}_{f}(v)=\varnothing$ for
all $v\in V_{f}^{++}$.
\end{enumerate}
\end{theorem}

We next summarize results on the above-mentioned parameters for spiders. Most
of the proofs can be found elsewhere, as indicated, and we only prove $(ii)$.

\begin{proposition}
\label{Prop_Spider}$(i)\hspace{0.1in}$\emph{\cite{MN}}\hspace{0.1in}For any
$r\geq1$ and $k\geq3$, $\alpha_{\operatorname{bn}}(\operatorname{Sp}%
(r^{k}))=kr$.

\begin{enumerate}
\item[$(ii)$] For any $r\geq2$ and $k\geq3$, $\alpha_{\operatorname{bnd}%
}(\operatorname{Sp}(r^{k}))=\alpha_{\operatorname{bnr}}(\operatorname{Sp}%
(r^{k}))=k(r-1)+1$.

\item[$(iii)$] \emph{\cite[Proposition 2.3.16]{LindaD}}\hspace{0.1in}For any
$k\geq3$, $\Gamma_{b}(\operatorname{Sp}(2^{k}))=k+1$.

\item[$(iv)$] \emph{\cite{Dunbar}\hspace{0.1in}}For any $r\geq1$ and $k\geq3$,
$\alpha_{h}(\operatorname{Sp}(r^{k}))=k(2r-1)$.
\end{enumerate}
\end{proposition}

\noindent\textbf{Proof of }$(ii)$\textbf{.}$\hspace{0.1in}$Let
$S=\operatorname{Sp}(r^{k})$. Denote the head of $S$ by $b$ and its legs by
$L_{1},...,L_{k}$. Say $L_{i}=b,v_{i,1},...,v_{i,r-1},\ell_{i}$, $i=1,...,k$.
Define the broadcast $f_{0}$ on $S$ by $f_{0}(\ell_{1})=r$, $f_{0}(\ell
_{i})=r-1$ for $2\leq i\leq k$, and $f_{0}(x)=0$ for all other vertices of
$S$. Then $b\in\operatorname{PB}_{f_{0}}(\ell_{1})$ while $b_{i}%
\in\operatorname{PB}_{f_{0}}(\ell_{i})$ for $2\leq i\leq k$, hence $f_{0}$ is
irredundant. It is clear that $f_{0}$ is a dominating broadcast; since $f_{0}$
is also irredundant, it is minimal dominating. Since no vertex is dominated by
more than one broadcasting vertex, $f_{0}$ is bn-independent. Therefore,
$f_{0}$ is a bnr-broadcast and a bnd-broadcast, from which it follows that
$\alpha_{\operatorname{bnr}}(S)\geq\alpha_{\operatorname{bnd}}(S)\geq
\sigma(f_{0})=k(r-1)+1$.

To prove the upper bound, let $F$ be the set of all $\alpha
_{\operatorname{bnr}}$-broadcasts $f^{\prime}$ such that no leaf hears
$f^{\prime}$ from a non-leaf vertex. By Theorem \ref{ThmTrees}$(ii)$, $F$ is
nonempty. Among all broadcasts in $F$, let $f$ be one having the fewest
non-leaf broadcasting vertices. 

Suppose $f$ has at least one non-leaf broadcasting vertex. Among all of these,
let $u$ be one nearest to a leaf. Since the proof works the same for $b$ and
for any vertex $v_{i,j}$ for some $i\in\{1,...,k\}$ and some $j\in
\{1,...,r-1\}$, we assume without loss of generality that $u=v_{1,j}$, where
we take $j=0$ if $u=b$. Since $v_{1,j}$ does not dominate $\ell_{1}$,
$f(\ell_{1})>0$. Say $f(\ell_{1})=t$. Since $\operatorname{PB}_{f}(\ell
_{1})\neq\varnothing$ and $v_{1,j}\in V_{f}^{+}$, the vertex $v_{1,r-t}$ on
$L_{1}$ belongs to $\operatorname{PB}_{f}(\ell_{1})$. By the choice of
$v_{1,j}$ as being a broadcasting vertex nearest to a leaf, either
$v_{1,r-t-1}$ is undominated, or $v_{1,r-t-1}\in\operatorname{PB}_{f}%
(v_{1,j})$. Since $v_{1,j}$ also does not dominate $\ell_{2}$, there exists a
vertex $x\in B_{f}(v_{1,j})$ that does not lie on the $\ell_{1}-v_{1,j}$
subpath of $L_{1}$; indeed, either (a) $x$ lies on $L_{1}$ between $v_{1,j}$
and $b$, or (b) $x=b$ if $f(v_{1,j})=d(v_{1,j},x)$, or (c) if $v_{1,j}$
overdominates $b$, then each leg $L_{i},\ 2\leq i\leq k$, contains such a
vertex $x$ (all at the same distance from $b$).

Let $d=d(v_{1,r-t},x)$ and note that $d\geq2f(v_{1,j})+1>f(v_{1,j})+1$. Define
the broadcast $g$ on $S$ by
\[
g(\ell_{1})=f(\ell_{1})+d-1,\ g(v_{1,j})=0,\ \text{and\ }%
g(v)=f(v)\ \text{otherwise.}%
\]
Then $\sigma(g)=\sigma(f)-f(v_{1,j})+d-1>\sigma(f)$. Moreover, since $\ell
_{1}$ does not broadcast to $x$, $N_{g}(\ell_{1})$ is a proper subset of
$N_{f}(\ell_{1})\cup N_{f}(v_{1,j})$, whereas $N_{g}(v)=N_{f}(v)$ for all
other vertices $v\in V_{g}^{+}$. This implies that $\operatorname{PB}%
_{g}(v)\neq\varnothing$ for each vertex in $V_{g}^{+}$; consequently, $g$ is a
bnr-broadcast. But $\sigma(g)>\sigma(f)=\alpha_{\operatorname{bnr}}(S)$, a
contradiction. We conclude that $V_{f}^{+}$ consists of leaves. 

Suppose some leaf $\ell$ overdominates $b$; assume without loss of generality
that $\ell=\ell_{1}$. Say $f(\ell_{1})=d(\ell_{1},b)+t=r+t$ for some $t$ such
that $1\leq t\leq r$. Since $f$ is irredundant and only leaves are
broadcasting vertices, it follows that $f(\ell_{i})\leq r-t-1$ for each
$i\in\{2,...,k\}$. Note that in this case, $\sigma(f)\leq
r+t+(k-1)(r-t-1)=k(r-1)-t(k-2)+1$. But then $\alpha_{\operatorname{bnr}%
}(S)=\sigma(f)<\sigma(f_{0})=k(r-1)+1$, where $f_{0}$ is the bnr-broadcast
defined in the first part of the proof. This is impossible. 

We deduce that no leaf overdominates $b$. Thus $f(\ell_{i})\leq r$ for each
$i$, and, since $f$ is irredundant and $\operatorname{PB}_{f}(\ell_{1}%
)\neq\varnothing$ for each $i$, $f(\ell_{i})=r$ for at most one $i$.
Therefore, $\sigma(f)\leq k(r-1)+1$. This proves that $\alpha
_{\operatorname{bnr}}(S)=k(r-1)+1$. Since $\alpha_{\operatorname{bnd}}%
(S)\leq\alpha_{\operatorname{bnr}}(S)$, the proof is complete.~$\blacksquare$

\bigskip

By Proposition \ref{Prop_Spider}, the differences $\alpha_{h}-\alpha
_{\operatorname{bn}}$, $\alpha_{h}-\alpha_{\operatorname{bnr}}$ and
$\alpha_{\operatorname{bn}}-\alpha_{\operatorname{bnr}}$ can be arbitrary,
because
\begin{align*}
\alpha_{h}(\operatorname{Sp}(r^{k}))-\alpha_{\operatorname{bn}}%
(\operatorname{Sp}(r^{k})) &  =k(2r-1)-kr=k\left(  r-1\right)  ,\\
\alpha_{h}(\operatorname{Sp}(r^{k}))-\alpha_{\operatorname{bnr}}%
(\operatorname{Sp}(r^{k})) &  =k(2r-1)-(kr-k+1)=kr-1\\
\text{and\ }\alpha_{\operatorname{bn}}(\operatorname{Sp}(r^{k}))-\alpha
_{\operatorname{bnr}}(\operatorname{Sp}(r^{k})) &  =kr-(kr-k+1)=k-1.
\end{align*}
Mynhardt and Neilson \cite{MN} proved that
\[
\alpha_{\operatorname{bn}}(G)/\alpha_{\operatorname{bnr}}(G)<2,\alpha
_{h}(G)/\alpha_{\operatorname{bn}}(G)<2\text{\ and\ }\alpha_{h}(G)/\alpha
_{\operatorname{bnr}}(G)<3
\]
for any graph $G$, and used spiders to illustrate that all bounds are
asymptotically best possible. We consider corresponding results for
$\alpha_{\operatorname{bn}}$ and $\alpha_{\operatorname{bnr}}$ versus
$\Gamma_{b}$ in Sections \ref{SecDiff} and \ref{SecRatios}.

\section{Non-comparability of parameters}

\label{Sec_Non_Comp}While it is clear from the definitions that $\alpha
_{\operatorname{bnd}}(G)\leq\Gamma_{b}(G)$ for all graphs $G$, we demonstrate
below that $\alpha_{\operatorname{bnr}}\diamond\Gamma_{b}$ and $\alpha
_{\operatorname{bn}}\diamond\Gamma_{b}$.%
\begin{figure}[ptb]%
\centering
\includegraphics[
height=1.3059in,
width=2.6662in
]%
{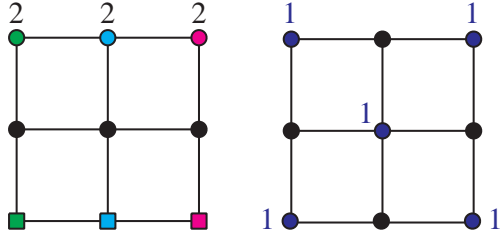}%
\caption{A $\Gamma_{b}$-broadcast of $G_{3,3}$ (left) of weight $6$ and an
$\alpha$-set (right) of cardinality $5$}%
\label{Fig3X3}%
\end{figure}

As shown in \cite{BF}, $\Gamma_{b}(G_{m,n})=m(n-1)$ for $2\leq m\leq n$. A
minimal dominating broadcast with weight $6$ of $G_{3,3}$ is depicted in
Figure \ref{Fig3X3}, where the square vertices belong to the private
boundaries of the broadcasting vertices of the same colour (on the same
vertical line when seen in monochrome). On the other hand, since grids are
$2$-connected bipartite graphs, (\ref{sequence5}) gives $\alpha(G_{3,3}%
)=\alpha_{\operatorname{bn}}(G_{3,3})=\alpha_{\operatorname{bnr}}%
(G_{3,3})=\left\lceil \frac{9}{2}\right\rceil =5<6=\Gamma_{b}(G_{3,3})$. For
the path $P_{n}$, $n\geq2$, $\alpha_{\operatorname{bnd}}(P_{n})=\alpha
_{\operatorname{bn}}(P_{n})=\alpha_{\operatorname{bnr}}(P_{n})=\Gamma
_{b}(P_{n})=n-1$, where the results for $\alpha_{\operatorname{bnd}}%
(P_{n}),\ \alpha_{\operatorname{bn}}(P_{n})$ and $\alpha_{\operatorname{bnr}%
}(P_{n})$ are stated in \cite[Section 3.2]{LindaD}, and the result for
$\Gamma_{b}(P_{n})$ follows from the bound of Bouchemakh and Fergani \cite{BF}.

\begin{figure}[ptb]%
\centering
\includegraphics[
height=2.4007in,
width=5.6879in
]%
{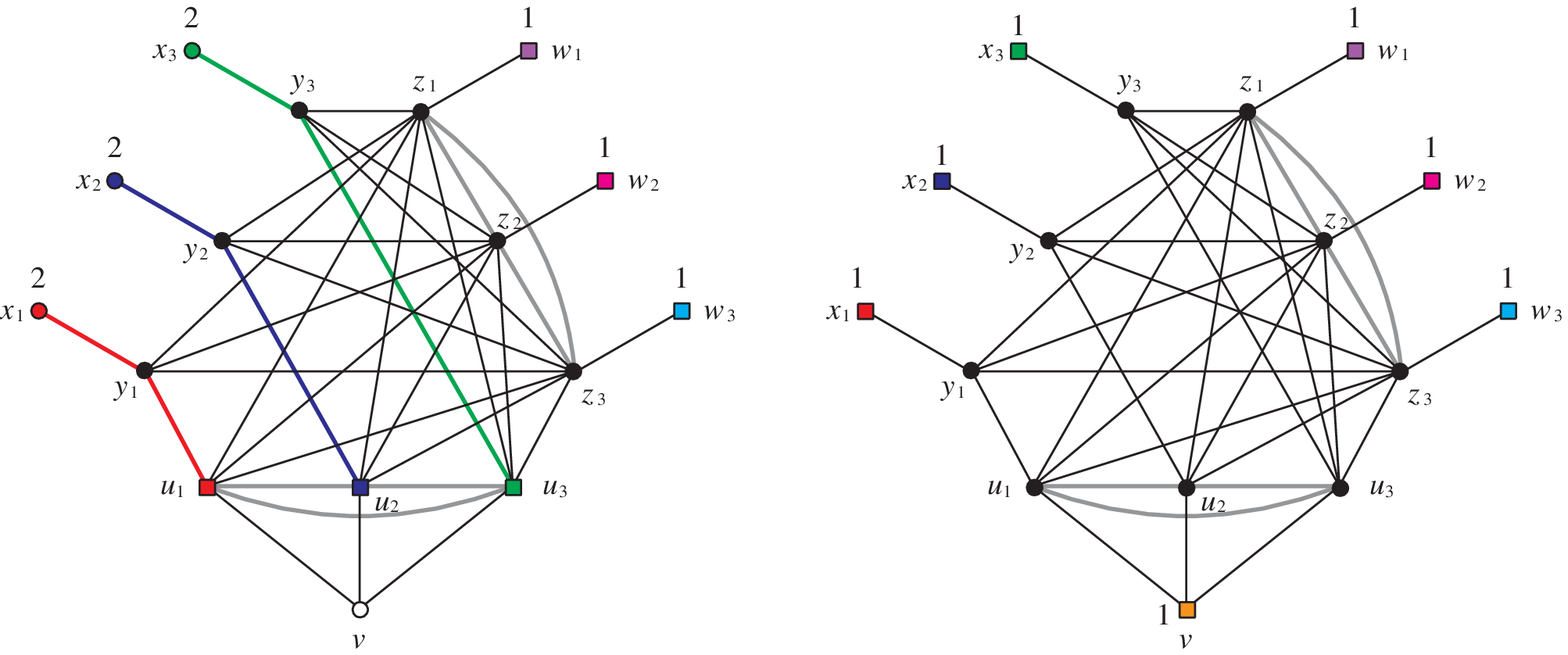}%
\caption{A graph $G$ with $\alpha_{\operatorname{bnr}}(G)=9$ and $\Gamma
_{b}(G)=7$. A non-dominating $\alpha_{\operatorname{bnr}}$-broadcast is shown
on the left, and a $\Gamma_{b}$-broadcast on the right}%
\label{Fig_alpha_bnr}%
\end{figure}

For an example of a graph for which $\alpha_{\operatorname{bn}}\geq
\alpha_{\operatorname{bnr}}>\Gamma_{b}$, consider the graph in Figure
\ref{Fig_alpha_bnr}. As shown in \cite{LindaD}, $\alpha_{\operatorname{bnr}%
}(G)\geq9$ (it can be shown that equality holds) and $\Gamma_{b}(G)=7$. We
omit the proofs here as they are special cases of the proofs of Propositions
\ref{bnr2b} and \ref{Gbnr} in Section \ref{Sec_bnr-Gam}. It follows that
$\alpha_{\operatorname{bnr}}\diamond\Gamma_{b}$ and $\alpha_{\operatorname{bn}%
}\diamond\Gamma_{b}$.%

\section{Differences}

\label{SecDiff}Having demonstrated that $\alpha_{\operatorname{bnr}}%
\diamond\Gamma_{b}$ and $\alpha_{\operatorname{bn}}\diamond\Gamma_{b}$, we
proceed to show that the differences $\Gamma_{b}-\alpha_{\operatorname{bn}}$
and $\Gamma_{b}-\alpha_{\operatorname{bnr}}$ can be arbitrary for general
graphs $G$, while $\Gamma_{b}-\alpha_{\operatorname{bnr}}$ can also be
arbitrary for trees. In the other direction we show that $\alpha
_{\operatorname{bnr}}-\Gamma_{b}$, hence also $\alpha_{\operatorname{bn}%
}-\Gamma_{b}$, can be arbitrary for trees as well as for cyclic graphs.

\subsection{$\Gamma_{b}(G)-\alpha_{\operatorname{bn}(\operatorname{bnr})}(G)$}

In this subsection we show that $\Gamma_{b}(G)-\alpha_{\operatorname{bn}}(G)$
is unbounded. Since $\alpha_{\operatorname{bn}}(G)\geq\alpha
_{\operatorname{bnr}}(G)$ it will follow that $\Gamma_{b}(G)-\alpha
_{\operatorname{bnr}}(G)$ is unbounded for graphs in general.

\begin{proposition}
\label{Prop_Gamma_minus_bn}{For any integer $k\geq3$, there exists a graph
$G_{k}$ such that $\Gamma_{b}(G_{k})-\alpha_{\operatorname{bn}}(G_{k}%
)\geq\lfloor\frac{k}{2}\rfloor$.}
\end{proposition}

\noindent\textbf{Proof.\hspace{0.1in}}Let $G_{k}$ be the grid $G_{3,k}=P_{3}%
\boksie
P_{k}$. Since grids are $2$-connected bipartite graphs, $\alpha
_{\operatorname{bn}}(G_{3,k})=\lceil\frac{3k}{2}\rceil$ by (\ref{sequence5}).
Let $X$ and $Y$ be the vertex sets of the fibres of $P_{k}$ that correspond to
leaves of $P_{3}$. Define the broadcast $f$ on $G_{k}$ by $f(x)=2$ when $x\in
X$ and $f(x)=0$ otherwise. For each $x\in X$, $\operatorname{PB}_{f}(x)$
consists of the vertex in $Y$ that belongs to the same fibre of $P_{3}$ as $x$
(see Figure \ref{Fig3X3} for $k=3$). Hence $f$ is an irredundant dominating
broadcast of weight $2k$. Thus $\Gamma_{b}(G_{3,k})\geq2k$ and $\Gamma
_{b}(G_{k})-\alpha_{\operatorname{bn}}(G_{k})\geq\lfloor\frac{k}{2}\rfloor
$.~$\blacksquare$

\bigskip

By Observation \ref{spanning_tree}, $\alpha_{\operatorname{bn}}(T)\geq
\alpha_{\operatorname{bn}}(G)$ when $T$ is a spanning tree of $G$. Hence it is
possible that $\Gamma_{b}(T)-\alpha_{\operatorname{bn}}(T)$ is bounded. We
leave this as an open question; see Question \ref{Q_bounded} in
Section~\ref{Sec_Qs}.

\begin{figure}[ptb]%
\centering
\includegraphics[
height=2.5322in,
width=5.8237in
]%
{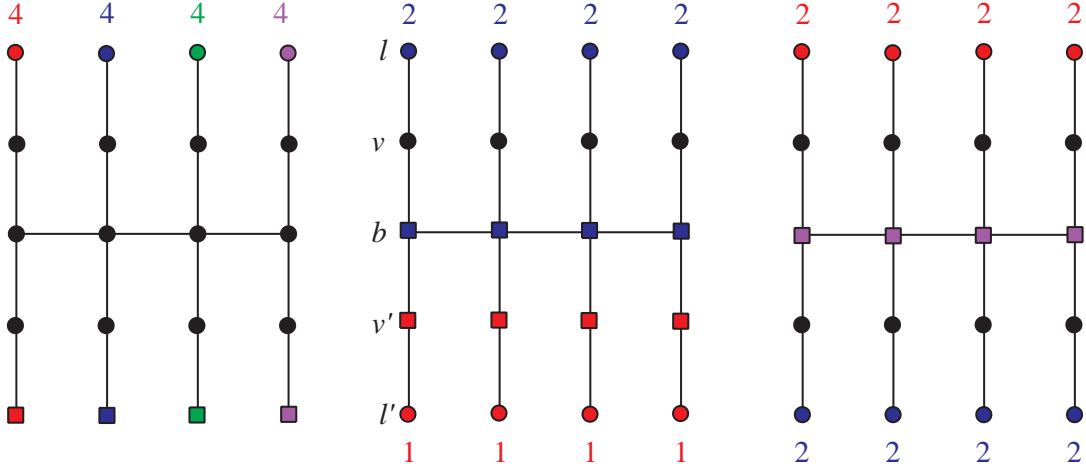}%
\caption{The tree $T_{4}$ with a minimal dominating broadcast on the left, an
{$\alpha_{\operatorname{bnr}}$-broadcast in the middle, and an }%
$\alpha_{\operatorname{bn}}$-broadcast on the right.}%
\label{Fig_Tk}%
\end{figure}

To show that $\Gamma_{b}(T)-\alpha_{\operatorname{bnr}}(T)$ is unbounded for
trees in general, we use the tree $T_{k}$ described below and shown in Figure
\ref{Fig_Tk} for $k=4$. To form $T_{k}$, take $k$ copies $P_{5}^{i}%
=(l_{i},v_{i},b_{i},v_{i}^{\prime},l_{i}^{\prime}),\ i=1,...,k,$ of $P_{5}$
and add the edges $b_{i}b_{i+1},\ i=1,...,k-1$. Define the broadcast $f$ on
$T_{k}$ by $f(l_{i})=4,\ i=1,...,k$, and $f(x)=0$ otherwise. Note that $f$ is
dominating and $\operatorname{PB}_{f}(l_{i})=\{l_{i}^{\prime}\}$. Therefore,
$f$ is minimal dominating, so $\Gamma_{b}(T_{k})\geq\sigma(f)=4k$.%

\begin{proposition}
{\label{bnrd} For $k\geq2$, $\alpha_{\operatorname{bnr}}(T_{k})=3k$, where
$T_{k}$ is the tree described above and shown in Figure \ref{Fig_Tk} (for
$k=4$). }
\end{proposition}

\noindent\textbf{Proof.\hspace{0.1in}}For $1\leq i\leq k$, let $A_{i}%
=V(P_{5}^{i})=\{l_{i},v_{i},b_{i},v_{i}^{\prime},l_{i}^{\prime}\}$. Define the
broadcast $g$ on $T_{k}$ by $g(l_{i})=2,\ i=1,...,k$, $g(l_{i}^{\prime
})=1,\ i=1,...,k$, and $g(x)=0$ otherwise. Then $g$ is bn-independent and
dominating. Furthermore, $\operatorname{PB}_{g}(l_{i})=\{b_{i}\}$ and
$\operatorname{PB}_{g}(l_{i}^{\prime})=\{l_{i}^{\prime},v_{i}^{\prime}\}$.
Therefore $g$ is bn-independent and irredundant, hence $\alpha
_{\operatorname{bnr}}(T_{k})\geq\sigma(g)=3k$.

We now show that $\alpha_{\operatorname{bnr}}(T_{k})\leq3k$. Suppose to the
contrary that $\alpha_{\operatorname{bnr}}(T_{k})>3k$. By Theorem
\ref{ThmTrees}$(ii)$, there exists an $\alpha_{\operatorname{bnr}}(T_{k}%
)$-broadcast such that no leaf hears a non-leaf. Of all such broadcasts, let
$f$ be one such that $L=\{v:v\in V(T_{k})\text{ and }f(v)\geq4\}$ is a
minimum. Let $M=\max\{f(v):v\in V_{f}^{+}\}$. Assume by symmetry that
$f(l_{i})=\max\{f(l_{i}),f(l_{i}^{\prime})\}$ for each $i$. Then
$f(l_{i}^{\prime})=\min\{f(l_{i}),f(l_{i}^{\prime})\}$. Since (a) no leaf
hears a non-leaf, (b) $f$ is bn-independent, and (c) $f$ is irredundant,
\begin{equation}%
\begin{tabular}
[c]{ll}%
(a) & $f(v_{i})=f(v_{i}^{\prime})=0$ and $f(b_{i})\leq1$ for $1\leq i\leq
k$,\\
(b) & $f(b_{i})+f(b_{i+1})\leq1$ for $1\leq i\leq k-1$,\\
(c) & $f(l_{i}^{\prime})\leq1$ for $1\leq i\leq k$.
\end{tabular}
\ \ \label{eq_props}%
\end{equation}

Moreover, the bn-independence of $f$ also implies that%
\begin{equation}
\text{if }f(b_{i})=1\text{ for some }i\text{, then }f(l_{i})=f(l_{i}^{\prime
})=1. \label{eq_prop2}%
\end{equation}

Suppose $L=\varnothing$. If $M<3$ then, by (\ref{eq_props}) and
(\ref{eq_prop2}), $f(A_{i})\leq3$ for all $1\leq i\leq k$. But then
$\sigma(f)\leq3k$, a contradiction. Hence assume that $M=3$. We obtain a
contradiction by showing that if $f(A_{i})>3$ for some $i$, then
$f(A_{1})+f(A_{2})\leq6$ if $i=1$, $f(A_{k-1})+f(A_{k})\leq6$ if $i=k$, and
$f(A_{i-1})+f(A_{i})+f(A_{i+1})\leq8$ otherwise.

By (\ref{eq_props}) and (\ref{eq_prop2}), if $1\leq f(l_{i})\leq2$ for some
$i$, then $f(A_{i})\leq3$. Assume therefore that $f(l_{i})=3$. If $i=1$, then
$B_{f}(l_{1})=\{v_{1}^{\prime},b_{2}\}$ and $f(b_{1})=f(b_{2})=0$. Now
$f(l_{1}^{\prime})\leq1$, otherwise $\operatorname{PB}_{f}(l_{1}^{\prime
})=\varnothing$. Similarly, $f(l_{2})\leq1$ and $f(l_{2}^{\prime})\leq1$.
Hence $f(A_{1})+f(A_{2})\leq6$. By symmetry, if $f(l_{k})=3$, then
$f(A_{k})+(A_{k-1})\leq6$. If $f(l_{i})=3$ and $i\neq1,k$, then $B_{f}%
(l_{i})=\{b_{i-1},v_{i}^{\prime},b_{i+1}\}$, $f(l_{i-1})=f(l_{i-1}^{\prime
})=f(l_{i}^{\prime})=f(l_{i+1})=f(l_{i+1}^{\prime})=1$ and $f(b_{i-1}%
)=f(b_{i})=f(b_{i+1})=0$. Hence $f(A_{i-1})+f(A_{i})+f(A_{i+1})=8$. It follows
that $\sigma(f)\leq3k$, contrary to our assumption.\label{March10}

Suppose $L\neq\varnothing$ and $M=4$.

\begin{itemize}
\item If $k=2$, then, without loss of generality, $f(l_{1})=4$. Then
$B_{f}(l_{1})=\{v_{2}^{\prime},v_{2},l_{1}^{\prime}\}$ and, since $f$ is
bnr-independent, $f(l_{2})=f(l_{2}^{\prime})=1$ and $f(b_{1})=f(b_{2}%
)=f(l_{1}^{\prime})=0$. Hence $f(A_{1})+f(A_{2})=6=3k$.

\item If $k=3$ and $f(l_{2})=4$, then $f(l_{1})=f(l_{1}^{\prime}%
)=f(l_{3})=f(l_{3}^{\prime})=1$ and $f(b_{1})=f(b_{2})=f(b_{3})=f(l_{2}%
^{\prime})=0$, hence $\sigma(f)=8<3k$.

\item If $f(l_{i})=4$ and either $k=3$ and $i\neq2$, or $k>3$, then, without
loss of generality, $\{l_{i}^{\prime},v_{i+1},v_{i+1}^{\prime},b_{i+2}%
\}\subseteq B_{f}(l_{i})$ and $f(l_{i+1})=f(l_{i+1}^{\prime})=1$ and
$f(b_{i})=f(b_{i+1})=f(b_{i+2})=0$. Create a new broadcast $g$ with
$g(l_{i})=3$, $g(l_{i}^{\prime})=1$ and $g(x)=f(x)$ otherwise. Notice that
$b_{i+1}\in\operatorname{PB}_{g}(l_{i})$. Since $g(l_{i}^{\prime})=1$,
$\operatorname{PB}_{g}(l_{i}^{\prime})=\{l_{i}^{\prime}\}$. Also, $N_{g}%
(l_{i})\cup N_{g}(l_{i}^{\prime})\subseteq N_{f}(l_{i})$. Hence, $g$ is a
bnr-independent broadcast. Since $\sigma(f)=\sigma(g)$, either $g$ can be
extended and violates the maximality of $f$, or, since it has fewer vertices
broadcasting with strength $4$, it violates the choice of $f$.
\end{itemize}

Assume that $M\geq5$. Let $l_{i}$ be a vertex such that $f(l_{i})=M$. Since
$f(l_{i})\leq e(l_{i})$ and by the structure of $T_{k}$, there are two leaves
$l_{t},l_{t}^{\prime}$, such that $d(l_{t},l_{i})=d(l_{t}^{\prime},l_{i})=M$.
Assume by symmetry that $t>i$. Since $d(l_{i},l_{i+1})=5$,
$t=(i+1)+(M-5)=i+M-4$. Create a new broadcast $g$ with $g(l_{j})=2$,
$g(l_{j}^{\prime})=1$ and $g(b_{j})=0$ for all $i\leq j\leq t$, and
$g(x)=f(x)$ otherwise. Notice that $\bigcup_{j=i}^{t}(N_{g}(l_{j})\cup
N_{g}(l_{j}^{\prime}))\subset N_{f}(l_{i})$. For $i\leq j\leq t$, $b_{j}%
\in\operatorname{PB}_{g}(l_{j})$ and $l_{j}^{\prime}\in\operatorname{PB}%
_{g}(l_{j}^{\prime})$. Hence $g$ is a bnr-independent broadcast with
\[
\sigma(g)=\sigma(f)-M+3(t-i+1)=\sigma(f)-M+3(M-3)=\sigma(f)+2M-9.
\]
Since $M\geq5$, $\sigma(g)>\sigma(f)$ and $g$ violates the maximality of $f$.
We conclude that $\alpha_{\operatorname{bnr}}(T_{k})\leq3k$.~$\blacksquare$

\bigskip

Since $\Gamma_{b}(T_{k})\geq4k$, the following result is an immediate
consequence of Proposition \ref{bnrd}.

\begin{theorem}
{\label{unb} For any integer $k\geq1$ there exists a tree $T$ such that
$\Gamma_{b}(T)-\alpha_{\operatorname{bnr}}(T)\geq k$. }
\end{theorem}

\subsection{$\alpha_{\operatorname{bnr}(\operatorname{bn})}(G)-\Gamma_{b}(G)$}

\label{Sec_bnr-Gam}

We next consider the differences $\alpha_{\operatorname{bn}}-\Gamma_{b}$ and
$\alpha_{\operatorname{bnr}}-\Gamma_{b}$ for both trees and cyclic graphs. By
Proposition \ref{Prop_Spider}, when $k\geq3$, then
\[
\alpha_{\operatorname{bn}}(\operatorname{Sp}(2^{k}))=2k\text{\ and\ }%
\alpha_{\operatorname{bnd}}(\operatorname{Sp}(2^{k}))=\alpha
_{\operatorname{bnr}}(\operatorname{Sp}(2^{k}))=\Gamma_{b}(\operatorname{Sp}%
(2^{k}))=k+1.
\]
Therefore $\alpha_{\operatorname{bn}}(\operatorname{Sp}(r^{k}))-\Gamma
_{b}(\operatorname{Sp}(r^{k}))=k-1$, and it follows that the difference
$\alpha_{\operatorname{bn}}-\Gamma_{b}$ can be arbitrary for trees. We
strengthen this result by constructing a tree $H_{k}$ such that $\alpha
_{\operatorname{bnr}}(H_{k})-\Gamma_{b}(H_{k})\geq k$.

\begin{figure}[ptb]%
\centering
\includegraphics[
height=4.4097in,
width=3.9583in
]%
{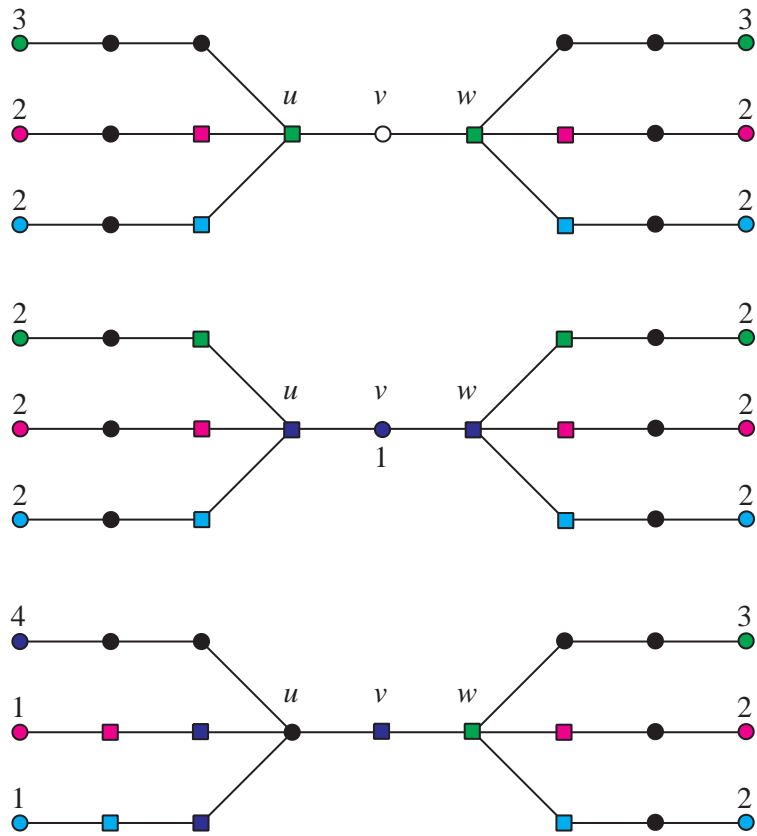}%
\caption{A tree $T$ with a non-dominating $\alpha_{\operatorname{bnr}}%
$-broadcast (top) and two $\Gamma_{b}$-broadcasts (middle and bottom)}%
\label{undombnr}%
\end{figure}

If a tree $T$ has an $\alpha_{\operatorname{bnr}}$-broadcast which is
dominating, then $\alpha_{\operatorname{bnr}}(T)\leq\Gamma_{b}(T)$. However,
not all trees have such a broadcast and there exist trees such that
$\alpha_{\operatorname{bnr}}(T)>\Gamma_{b}(T)$. Figure \ref{undombnr} gives an
example of a bnr-independent broadcast on a tree $T$ which is not dominating.
By using symmetry and examining a few cases, it can be shown that
$\alpha_{\operatorname{bnr}}(T)=14$ and $\Gamma_{b}(T)=13$. We use this tree
as a basis to construct bigger trees $H_{k}$ mentioned above.

To construct $H_{k}$, take $3k$ copies $T_{1},...,T_{3k}$ of the tree in
Figure \ref{undombnr} and label the central vertex and its neighbours in the
$i^{\operatorname{th}}$ copy as $u_{i},v_{i},w_{i}$. Let $H_{k}$ be the tree
formed by joining $v_{i}$ to $v_{i+1}$ for each $i=1,...,3k-1$. For
$i=1,...,3k$, let $f_{i}$ be the broadcast on $T_{i}$ illustrated in the top
copy of $T$ in Figure \ref{undombnr}, and define $f=\bigcup_{i=1}^{3k}f_{i}$.
Then $f$ is a bnr-broadcast, hence $\alpha_{\operatorname{bnr}}(H_{k})\geq
42k$. We show that ${\Gamma_{b}}${$(H_{k})=41k$.}

\begin{proposition}
\label{Prop_Hk}{For $k\geq1$, }${\Gamma_{b}}${$(H_{k})=41k$, where $H_{k}$ is
the tree described above.}
\end{proposition}

\noindent\textbf{Proof.\hspace{0.1in}}Let $f_{i},\ g_{i}$ and $h_{i}$ be the
broadcasts on $T_{i}$ illustrated in the top, middle and bottom copy,
respectively, of $T$ in Figure \ref{undombnr}. Define the broadcast $\lambda$
on $H_{k}$ by%
\[
\lambda(x)=\left\{
\begin{tabular}
[c]{ll}%
$g_{i}(x)$ & if $x\in V(T_{i})$ and $i\equiv2\ (\operatorname{mod}\ 3)$\\
$f_{i}(x)$ & otherwise.
\end{tabular}
\ \ \ \ \right.
\]
Since $\sigma(f_{i})=14$ and $\sigma(g_{i})=13$, $\sigma(\lambda
)=28k+13k=41k$. It is easy to see that $\lambda$ is a dominating broadcast.
Suppose $i\equiv2\ (\operatorname{mod}\ 3)$. Then $\operatorname{PB}_{\lambda
}(v_{i})=\{u_{i},v_{i},w_{i},v_{i-1},v_{i+1}\}$, and if $\ell$ is a leaf, then
$\operatorname{PB}_{\lambda}(\ell)$ consists of the vertex at distance $2$
from $\ell$. Suppose $i\equiv0$ or $1\ (\operatorname{mod}\ 3)$. If
$\lambda(\ell)=3$, then $\operatorname{PB}_{g}(\ell)=\{u_{i}\}$ or $\{w_{i}%
\}$, as the case may be, and if $\lambda(\ell)=2$, then, as before,
$\operatorname{PB}_{\lambda}(\ell)$ consists of the vertex at distance $2$
from $\ell$. Hence $\lambda$ is a minimal dominating broadcast and we deduce
that $\Gamma_{b}(H_{k})\geq41k$. We show that $\Gamma_{b}(H_{k})=41k$.%

Among all $\Gamma_{b}(H_{k})$-broadcasts on $H_{k}$, let $\rho$ be one such
that $V_{\rho}^{+}$ contains the maximum number of leaves. For $i=1,...,3k$,
define the function $\rho_{i}:V(T_{i})\rightarrow\{1,...,\operatorname{diam}%
(H_{k})\}$ by $\rho_{i}(x)=\rho(x)$ for each $x\in V(T_{i})$. (Possibly,
$\rho_{i}$ is not, strictly speaking, a broadcast on $T_{i}$ because $\rho
_{i}(x)$ may exceed the eccentricity of $x$ in $T_{i}$.) If each $v_{i}$ is
dominated only by a vertex of $T_{i}$, then $\rho_{i}$ is a broadcast on
$T_{i}$ and, as in the case of $T$, $\sigma(\rho_{i})\leq13$ for each $i$, so
that $\sigma(\rho)\leq39k<41k$. Hence assume a vertex $v_{i}$ is dominated by
a vertex $x\in V(T_{j})$, $j\neq i$.

\begin{claim}
\label{Cl1}If $x\in V(T_{j})$ dominates $v_{i}$, where $j\neq i$, then $x$
does not overdominate~$v_{i}$.
\end{claim}

\noindent\textbf{Proof of Claim \ref{Cl1}\hspace{0.1in}}Say $x$ overdominates
$v_{i}$ by exactly $t$, where$\ t\geq0$, that is, $\rho(x)=t+d(x,v_{i})$. The
possible values of $\sigma(\rho_{i})$ and $\sigma(\rho_{j})$ are tabled below.%
\[%
\begin{tabular}
[c]{|c|l|l|l|l|l|}\hline
$t$ & $0$ & $1$ & $2$ & $3$ & $\geq4$\\\hline
\multicolumn{1}{|l|}{$\sigma(\rho_{i})$} & $\leq14$ & $\leq12$ & $\leq6$ &
$\leq6$ & $0$\\\hline
\multicolumn{1}{|l|}{$\sigma(\rho_{j})$} & $\leq12+\rho(x)$ & $\leq6+\rho(x)$
& $\leq6+\rho(x)$ & $\rho(x)$ & $\rho(x)$\\\hline
\end{tabular}
\ \ \
\]

Suppose $t>0$. Let $j$ be the smallest index such that $T_{j}$ contains a
vertex $x$ which overdominates a vertex $v_{l},\ l\neq j$. By symmetry we may
assume that $j\leq\left\lceil \frac{3k}{2}\right\rceil $. We show that $x$
overdominates some vertex $v_{i}$ by exactly $1$.

\begin{itemize}
\item Suppose $x$ overdominates $v_{3k}$ by at least $4$. Since $j\leq
\left\lceil \frac{3k}{2}\right\rceil $, $x$ overdominates each $v_{i}%
,\ i=1,...,3k$, by at least $4$. Then, regardless of the value of $j$, $x$
dominates $H_{k}$ and $\sigma(\rho)=e(x)\leq\operatorname{diam}(H_{k}%
)=3k+7<41k$.

\item Suppose $x$ overdominates $v_{3k}$ by $3$.

\begin{itemize}
\item If $k$ is odd and $j=\left\lceil \frac{3k}{2}\right\rceil $, then
$\rho(x)\leq\left\lfloor \frac{3k}{2}\right\rfloor +7$ and $x$ dominates all
of $H_{k}$ except for the twelve leaves of $T_{1}$ and $T_{3k}$, so
$\sigma(\rho)\leq\left\lfloor \frac{3k}{2}\right\rfloor +19<41k$ for all $k$.

\item In all other cases, $\rho(x)\leq\operatorname{diam}(H_{k})-1$ and $x$
dominates all of $H_{k}$ except for six leaves, so $\sigma(\rho)\leq3k+12<41k$
for all $k$.
\end{itemize}

\item Similarly, if $x$ overdominates $v_{3k}$ by $2$, we obtain that
$\sigma(\rho)\leq\operatorname{diam}(H_{k})-2+24=3k+29<41k$ for all $k$.
\end{itemize}

In each case, $\sigma(\rho)<\sigma(\lambda)$ and we have a contradiction of
the assumption that $\rho$ is a $\Gamma_{b}(H_{k})$-broadcast. We deduce that
$x$ does not overdominate $v_{3k}$ by more than $1$, that is, $\rho(x)\leq
d(x,v_{3k})+1$. Say $\rho(x)=d(x,v_{3k})+1-\beta$, for some $\beta$ such that
$0\leq\beta\leq3k-j-1$. That is, $\rho(x)=3k-j-\beta+1$. Observe that if
$\beta=0$, then $x$ overdominates $v_{3k}$ by exactly $1$, if $\beta=1$, then
$x$ overdominates $v_{3k-1}$ by exactly $1$, etc., and if $\beta=3k-j-1$, then
$x$ overdominates $v_{j+1}$ by exactly $1$. By symmetry and the choice of $j$,
it is also possible -- but not necessary -- that $x$ overdominates some vertex
$v_{i^{\prime}}$ by exactly $1$, where $i^{\prime}<j$. Since $x$ overdominates
at least one $v_{i}$ by exactly $1$, there exists a smallest index $i\neq j$
such that $x$ overdominates $v_{i}$ by exactly $1$. We consider two cases,
depending on the value of $i$ relative to $j$.\smallskip

\noindent\textbf{Case 1:\hspace{0.1in}}$i<j$. Since $j\leq\left\lceil
\frac{3k}{2}\right\rceil $, by symmetry $x$ also overdominates $v_{2j-i}$ by
exactly $1$. Then $\sigma(\rho_{i}),\sigma(\rho_{2j-i})\leq12$. If the indices
$i-1$ and $2j-i+1$ exist, $x$ dominates $v_{i-1}$ and $v_{2j-i+1}$. In this
case, $\sigma(\rho_{i-1}),\sigma(\rho_{2j-i+1})\leq14$. It is possible that
one or both of $v_{i-1}$ and $v_{2j-i+1}$ hear $\rho$ from vertices different
from $x$. Define the broadcast $\eta$ on $H_{k}$ by%
\[
\eta(y)=\left\{
\begin{tabular}
[c]{ll}%
$h_{\lambda}(y)$ & if $y\in V(T_{\lambda})$ and $\lambda\in\{i-1,2j-i+1\}$\\
$f_{\lambda}(y)$ & if $y\in V(T_{\lambda})$ and $\lambda\in\{i,2j-i\}$\\
$\rho(y)-1$ & if $y=x$\\
$\rho(y)$ & otherwise.
\end{tabular}
\ \right.
\]
Then
\begin{align*}
\sigma(\eta) &  =\sigma(\rho)-[\sigma(\rho_{i})+\sigma(\rho_{2j-i}%
)+\sigma(\rho_{i-1})+\sigma(\rho_{2j-i+1})]\\
&  +[\sigma(f_{i})+\sigma(f_{2j-i})+\sigma(h_{i-1})+\sigma(h_{2j-i+1})]-1.
\end{align*}
Since $\sigma(f_{i})-\sigma(\rho_{i})\geq2$ and $\sigma(\rho_{i-1}%
)-\sigma(h_{i-1})\leq1$ (and similarly for the other indices), we see that
$\sigma(\eta)>\sigma(\rho)$. Note that $v_{i}\in\operatorname{PB}_{\eta}(x)$.
The private boundaries of the vertices in $V(T_{i-1})\cap V_{\eta}^{+}$ and
$V(T_{2j-i+1})\cap V_{\eta}^{+}$ are as shown in Figure \ref{undombnr} but
possibly excluding $v_{i-1}$ or $v_{2j-i+1}$. Hence $\eta$ is bnr-independent.
All vertices $v_{l}$ are dominated, so all undominated vertices lie on
endpaths. Since $x$ $\rho$-dominates these vertices and $v_{i}\in
\operatorname{PB}_{\eta}(x)$, these vertices can be dominated by extending
$\eta$ on an appropriate leaf to get a minimal dominating broadcast
$\eta^{\ast}$ such that $\sigma(\rho)<\sigma(\eta^{\ast})$. But again,
this$\ $contradicts $\rho$ being a $\Gamma_{b}(H_{k})$-broadcast.\smallskip

\noindent\textbf{Case 2:\hspace{0.1in}}$i>j$. By the choice of $i$, $x$
overdominates each $v_{l}$, $1\leq l<i$, by at least $2$. Now $\sigma(\rho
_{i})\leq12$ and $\sigma(\rho_{i-1})=6$; moreover, if $i<3k$, then
$\sigma(\rho_{i+1})\leq14$. Consider an endpath $P$ in $T_{i-1}$ from
$u_{i-1}$ to a leaf $\ell$, say $P=(u_{i-1},a,c,\ell)$. Then $a\in B_{\rho
}(x)$. If $a\notin\operatorname{PB}_{\rho}(x)$, then $\rho(c)=1$ and thus
$\rho(\ell)=0$. But then $\rho^{\prime}=(\rho-\{(c,1),(\ell,0)\})\cup
\{\{(c,0),(\ell,1)\}$ is a $\Gamma_{b}(H_{k})$-broadcast such that
$V_{\rho^{\prime}}^{+}$ contains more leaves than $V_{\rho}^{+}$ does,
contrary to the choice of $\rho$. We deduce that $a\in\operatorname{PB}_{\rho
}(x)$ and $\rho(\ell)=1$. Define the broadcast $\eta^{\prime}$ on $H_{k}$ by%
\[
\eta^{\prime}(y)=\left\{
\begin{tabular}
[c]{ll}%
$h_{i+1}(y)$ & if $y\in V(T_{i+1})$\\
$f_{i}(y)$ & if $y\in V(T_{i})$\\
$\rho(y)-1$ & if $y=x$\\
$\rho(y)$ & otherwise.
\end{tabular}
\ \ \right.
\]
Then
\[
\sigma(\eta^{\prime})=\sigma(\rho)-[\sigma(\rho_{i})+\sigma(\rho
_{i+1})]+[\sigma(f_{i})+\sigma(h_{i+1})]-1.
\]
As above, $\sigma(f_{i})-\sigma(\rho_{i})\geq2$ and $\sigma(\rho_{i+1}%
)-\sigma(h_{i+1})\leq1$, hence $\sigma(\rho)\leq\sigma(\eta^{\prime})$.
Moreover, since $a\in\operatorname{PB}_{\rho}(x)$, $\eta^{\prime}$ is not
dominating. Following the reasoning above we can extend $\eta^{\prime}$ on
appropriate leaves (including $\ell$) to get a minimal dominating broadcast
$\eta^{\prime\ast}$ such that $\sigma(\rho)<\sigma(\eta^{\prime\ast})$, again
a contradiction.~$\blacklozenge$

\bigskip

Therefore, if $x\in V(T_{j})$ dominates $v_{i}$ for $i\neq j$, then $i=j\pm1$
and $\{v_{j-1},v_{j+1}\}\subseteq B_{\rho}(x)$. (Assume $j-1$ and $j+1$ both
exist; the proof is the same if only one of them exists.) Then $\sigma
(\rho_{j})\leq13$ and $\sigma(\rho_{j-1}),\sigma(\rho_{j+1})\leq14$.
Consequently,
\[
\sigma(\rho_{l})\leq\left\{
\begin{tabular}
[c]{ll}%
$13$ & if some vertex of $T_{l}$ dominates a vertex of $T_{i},\ i\neq l$\\
$14$ & otherwise.
\end{tabular}
\ \ \ \ \right.
\]
Since each $v_{i}$ is dominated and the subtree of $T$ induced by
$\{v_{1},...,v_{3k}\}$ is the path $P_{3k}$, there are at least $\gamma
(P_{3k})=k$ indices $j$ such that some vertex of $T_{j}$ dominates
$v_{i},\ i\neq j$. This implies that $\sigma(\rho)\leq13k+14\cdot2k=41k$ and
we conclude that $\Gamma_{b}(H_{k})=41k$.~$\blacksquare$

\bigskip

Since $\alpha_{\operatorname{bnr}}(H_{k})\geq42k$ and $\Gamma_{b}(H_{k})=41k$,
the next theorem follows.

\begin{theorem}
For any integer $k\geq1$ there exists a tree $H_{k}$ such that $\alpha
_{\operatorname{bn}}(H_{k})-\Gamma_{b}(H_{k})\geq\alpha_{\operatorname{bnr}%
}(H_{k})-\Gamma_{b}(H_{k})\geq k$.
\end{theorem}

By Observation \ref{spanning_tree}, $\alpha_{\operatorname{bn}}(G)\leq
\alpha_{\operatorname{bn}}(T)$ if $T$ is a spanning tree of $G$. Therefore it
is possible that $\alpha_{\operatorname{bn}}(G)-\Gamma_{b}(G)$ is bounded for
cyclic graphs. Again, we use the unboundedness of $\alpha_{\operatorname{bnr}%
}(G)-\Gamma_{b}(G)$ to show that this is not the case. To show that
$\alpha_{\operatorname{bnr}}(G)-\Gamma_{b}(G)$ is unbounded, we generalize the
construction of the graph $G$ in Figure \ref{Fig_alpha_bnr}.

Denote the corona of $G$ and $K_{1}$ by $G\circ K_{1}$. For $k\geq1$,
construct the graph $G_{k}$ as follows. Let $U=\{u_{1},...,u_{k+1}\}$,
$W=\{w_{1},...,w_{k+1}\}$, $X=\{x_{1},...,x_{k+1}\}$, $Y=\{y_{1}%
,...,y_{k+1}\}$, $Z=\{z_{1},...,z_{k+1}\}$ and $\{v\}$ be disjoint sets of
vertices. Add edges so that
\begin{align*}
G_{k}[X]  &  \cong G_{k}[Y]\cong G_{k}[W]\cong\overline{K_{k+1}},\\
G_{k}[U\cup\{v\}]  &  \cong K_{k+2},\ G_{k}[Z]\cong K_{k+1},\ G_{k}[U\cup
Z]\cong K_{2(k+1)},\\
G_{k}[\{y_{i}\}\cup Z]  &  \cong K_{k+2}\ \text{for each }i\in\{1,...,k+1\},\\
G_{k}[Y\cup U]  &  \cong G_{k}[W\cup Z]\cong K_{k+1}\circ K_{1},\\
G_{k}[X\cup Y]  &  \cong\overline{K_{k+1}}%
\boksie
K_{2}\cong(k+1)K_{2}.
\end{align*}
Assume that the perfect matchings of $G_{k}[U\cup Y]$, $G_{k}[X\cup Y]$ and
$G_{k}[W\cup Z]$ are $\{u_{i}y_{i}:i=1,...,k+1\}$, $\{x_{i}y_{i}%
:i=1,...,k+1\}$ and $\{w_{i}z_{i}:i=1,...,k+1\}$, respectively. The graph
$G_{2}$ is illustrated in Figure \ref{Fig_alpha_bnr}.

\begin{proposition}
\label{bnr2b} Let $G_{k}$ be the graph described above and shown in Figure
\ref{Fig_alpha_bnr} for $k=2$. For any integer $k\geq1$, $\alpha
_{\operatorname{bnr}}(G_{k})\geq3(k+1)$.
\end{proposition}

\noindent\textbf{Proof.\hspace{0.1in}}Define the broadcast $f$ by%
\[
f(x)=\left\{
\begin{tabular}
[c]{ll}%
$2$ & if $x\in X$\\
$1$ & if $x\in W$\\
$0$ & otherwise.
\end{tabular}
\ \ \ \ \right.
\]
Then $\sigma(f)=3(k+1)$ and%
\[%
\begin{tabular}
[c]{rll}%
$N_{f}(x_{i})$ & $=$ & $\{x_{i},y_{i},u_{i}\}\cup Z$\\
$B_{f}(x_{i})$ & $=$ & $\{u_{i}\}\cup Z$\\
$\operatorname{PB}_{f}(x_{i})$ & $=$ & $\{u_{i}\}$\\
$N_{f}(w_{i})$ & $=$ & $\{w_{i},z_{i}\}$\\
$\operatorname{PB}_{f}(w_{i})$ & $=$ & $\{w_{i}\}$%
\end{tabular}
\]
for each $i\in\{1,...,k+1$. Thus we see that $f$ is a bnr-independent
broadcast, hence $\alpha_{\operatorname{bnr}}(G_{k})\geq3(k+1)$.~$\blacksquare
$

\begin{proposition}
{\label{Gbnr}} Let $G_{k}$ be the graph described above and shown in Figure
\ref{Fig_alpha_bnr} for $k=2$. For any integer $k\geq1${, $\Gamma_{b}%
(G_{k})=2k+3$.}
\end{proposition}

\noindent\textbf{Proof.\hspace{0.1in}}The set $X\cup W\cup\{v\}$ is an
independent dominating set of $G_{k}$ of cardinality $2k+3$, and its
characteristic function is a minimal dominating broadcast. Hence {$\Gamma
_{b}(G_{k})\geq2k+3$. }

Consider any minimal dominating broadcast $f$ on $G_{k}$. By symmetry, there
are six possible ways to dominate the vertex $v$.

\noindent\textbf{Case 1:}\hspace{0.1in}$v$ dominates itself. If $f(v)=3=e(v)$,
then $N_{f}(v)=V(G_{k})$, hence $\sigma(f)=3$. Suppose $f(v)=2$. Then
$N_{f}(v)=V(G_{k})-X-W$. Since $f$ is irredundant there exists an index $l$
such that $\{y_{l},z_{l}\}\cap\operatorname{PB}_{f}(v)\neq\varnothing$. Hence
either $f(Z)=f(y_{l})=f(x_{l})=f(u_{l})=0$ (so that $y_{l}\in\operatorname{PB}%
_{f}(v)$) or $f(U\cup Y)=f(z_{l})=f(w_{l})=0$ (so that $z_{l}\in
\operatorname{PB}_{f}(v)$) or both. But then $x_{l}$ or $w_{l}$ is not
$f$-dominated, a contradiction. If $f(v)=1$, then $f(U)=0$. In addition,
$f(z_{j})+f(w_{j})\leq1$ and $f(x_{j})+f(y_{j})\leq1$ for all $1\leq j\leq
k+1$. Hence $\sigma(f)\leq2(k+1)+1=2k+3$.\smallskip

\noindent\textbf{Case 2:}\hspace{0.1in}$v$ is dominated by $u_{l}$ for some
$l$. If $f(u_{l})=3=e(u_{l})$, then $u_{l}$ dominates $G_{k}$ and
$\sigma(f)=3$. Suppose $f(u_{l})=2$. Then $N_{f}(u_{l})=(V(G_{k}%
)-X)\cup\{x_{l}\}$. Since $f$ is irredundant, $f(v)=0$ and $f(z_{i}%
)=f(w_{i})=0$ for all $1\leq i\leq k+1$. Moreover, $f(u_{i})=0$ or
$f(u_{i})=2$ for $i\neq l$. Let $U=\{u_{i}:f(u_{i})=2\}$. For all $1\leq j\leq
k+1$, if $u_{j}\in U$ then $f(x_{j})=0$, and if $u_{j}\notin U$ then
$f(x_{j})+f(y_{j})\leq1$. Hence $\sigma(f)\leq2|U|+k+1-|U|=|U|+k+1\leq2k+2$.
Suppose $f(u_{l})=1$. Then $N_{f}(u_{l})=U\cup Z\cup\{v,y_{l}\}$. If
$f(u_{l^{\prime}})=1$ for $l^{\prime}\neq l$, then $\operatorname{PB}%
_{f}(u_{l})=\{y_{l}\}$ and $x_{l}$ is not $f$-dominated, a contradiction.
Hence $f(U)=1$. Moreover, $f(z_{j})+f(w_{j})\leq1$ and $f(x_{j})+f(y_{j}%
)\leq1$ for all $1\leq j\leq k+1$. Hence $\sigma(f)\leq2(k+1)+1=2k+3$%
.\smallskip

\noindent\textbf{Cases 3--4:}\hspace{0.1in}$v$ is dominated by $s\in
\{x_{l},y_{l}\}$ for some $l$, $1\leq l\leq k+1$. In either case
$N_{f}(s)=(V(G_{k})-X)\cup\{x_{l}\}$. By irredundance, $f(x_{i})\leq1$ for
$i\neq l$, $f(s)\leq e(s)$ and $f(x)=0$ otherwise. Hence $\sigma(f)\leq
f(s)+k\leq4+k\leq2k+3$.\smallskip

\noindent\textbf{Cases 5--6:}\hspace{0.1in}$v$ is dominated by $s\in
\{z_{l},w_{l}\}$ for some $l$, $1\leq l\leq k+1$. Since $d(s,v)=e(s)$,
$N_{f}(s)=V(G_{k})$. Hence $\sigma(f)=f(s)=e(s)\leq3$.\smallskip

This exhausts all possibilities, hence $\Gamma_{b}(G_{k})\leq2k+3$%
.~$\blacksquare$

\bigskip

Propositions \ref{bnr2b} and \ref{Gbnr} imply the following theorem.

\begin{theorem}
\label{Thm_cyclic}{For any integer $k\geq1$ there exists a cyclic graph
$G_{k}$ such that }$\alpha_{\operatorname{bn}}(G_{k})-\Gamma_{b}(G_{k}%
)\geq\alpha${$_{\operatorname{bnr}}(G_{k})-\Gamma_{b}(G_{k})\geq k$. }
\end{theorem}

\section{Ratios}

\label{SecRatios}

\subsection{$\alpha_{\operatorname{bnr}(\operatorname{bn})}(G)/\Gamma_{b}(G)$}

We show that the ratios $\alpha_{\operatorname{bn}}(G)/\Gamma_{b}(G)$ and
$\alpha_{\operatorname{bnr}}(G)/\Gamma_{b}(G)$ are bounded. We need the
following lemma.

\begin{lemma}
\label{Lem_reduce}Let $f$ be an $\alpha_{\operatorname{bn}}$-broadcast on a
graph $G$ such that $\operatorname{PB}_{f}(v)=\varnothing$ for some $v\in
V_{f}^{+}$. The broadcast $g$ on $G$ defined by $g(v)=f(v)-1$%
\ and\ $g(x)=f(x)$\ otherwise, is a dominating bn-independent broadcast, $v\in
V_{g}^{+}$ and $B_{g}(v)=\operatorname{PB}_{g}(v)\neq\varnothing$.
\end{lemma}

\noindent\textbf{Proof.\hspace{0.1in}}By Observation \ref{Ob_PB}, $v\in
V_{f}^{++}$, hence $v\in V_{g}^{+}$. If some vertex $u$ of $G$ is
$g$-undominated, then $u\in\operatorname{PB}_{f}(v)$, a contradiction because
$\operatorname{PB}_{f}(v)=\varnothing$. Hence $g$ is dominating. Since
$f(v)\leq e(v)$, $B_{f}(v)\neq\varnothing$. Therefore $B_{g}(v)\neq
\varnothing$. Consider any $u\in B_{g}(v)$ and suppose $u$ also hears $g$ from
$w\in V_{g}^{+}-\{v\}$. Let $x$ be any neighbour of $u$ on a $u-w$ geodesic.
Then $x$ hears $f$ from $v$, which means that $ux$ is $f$-covered by both $v$
and $w$, which is impossible because $f$ is bn-independent. Therefore
$u\in\operatorname{PB}_{g}(v)$. It follows that $B_{g}(v)=\operatorname{PB}%
_{g}(v)\neq\varnothing$.~$\blacksquare$

\begin{theorem}
\label{gammabound}{For any graph $G$,
\[
\alpha_{\operatorname{bnr}}(G)/\Gamma_{b}(G)\leq\alpha_{\operatorname{bn}%
}(G)/\Gamma_{b}(G)<2.
\]
The bound for }$\alpha_{\operatorname{bn}}(G)/\Gamma_{b}(G)$ is asymptotically
best possible.
\end{theorem}

\noindent\textbf{Proof.\hspace{0.1in}}If there exists an $\alpha
_{\operatorname{bn}}(G)$-broadcast that is minimal dominating, then
$\alpha_{\operatorname{bn}}(G)\leq\Gamma_{b}$. Hence we assume that no such
broadcast exists. Every $\alpha_{\operatorname{bn}}$-broadcast is dominating,
and an irredundant dominating broadcast is minimal dominating. Thus we may
also assume that no $\alpha_{\operatorname{bn}}$-broadcast on $G$ is
irredundant. We describe a strategy for turning a non-irredundant dominating
$\alpha_{\operatorname{bn}}$-broadcast into an irredundant dominating
broadcast with weight large enough to achieve the desired result.

Consider an $\alpha_{\operatorname{bn}}$-broadcast $f$ on $G$. Since $f$ is
not irredundant, there exists a vertex $v_{1}\in V_{f}^{+}$ such that
$\operatorname{PB}_{f}(v_{1})=\varnothing$. By Observation \ref{Ob_PB},
$v_{1}\in V_{f}^{++}$ and $|V_{f}^{+}|\geq2$. Define the broadcast $f_{1}$ by
\[
f_{1}(v_{1})=f(v_{1})-1\text{\ and\ }f_{1}(x)=f(x)\text{\ otherwise.}%
\]
By Lemma \ref{Lem_reduce}, $f_{1}$ is a dominating bn-independent broadcast.
Notice that $V_{f_{1}}^{+}=V_{f}^{+}$. Since $f(v_{1})\leq e(v_{1})$,
$B_{f}(v_{1})\neq\varnothing$. Therefore the set $B_{f_{1}}(v_{1}%
)=\{v:d(v,v_{1})=f(v_{1})-1\}$ is also nonempty. By the bn-independence of $f$
and the definition of $f_{1}$, $\operatorname{PB}_{f_{1}}(v_{1})=B_{f_{1}%
}(v_{1})\neq\varnothing$. Moreover, for each $u\in V_{f}^{+}-\{v_{1}\}$,
$\operatorname{PB}_{f}(u)\subseteq\operatorname{PB}_{f_{1}}(u)$, so if
$\operatorname{PB}_{f}(u)\neq\varnothing$, then $\operatorname{PB}_{f_{1}%
}(u)\neq\varnothing$. We conclude that $V_{f_{1}}^{+}$ contains more vertices
with non-empty private boundaries than $V_{f}^{+}$ does. 

If $f_{1}$ is not irredundant, we repeat this process, choosing a vertex
$v_{2}\in V_{f_{1}}^{++}-\{v_{1}\}$ having $\operatorname{PB}_{f_{1}}%
(v_{2})=\varnothing$, until we have a smallest index $k\geq1$ such that
$f_{k}$ is a dominating irredundant broadcast. Then $\sigma(f_{k})\leq
\Gamma_{b}(G)$. We show that $\sigma(f_{k})>\frac{1}{2}\sigma(f)$. 

Clearly, $\sigma(f)\geq|V_{f}^{1}|+2|V_{f}^{++}|$, hence $|V_{f}^{++}%
|\leq\frac{1}{2}\sigma(f)$. 

\begin{itemize}
\item If $V_{f}^{1}\neq\varnothing$, then $|V_{f}^{++}|<\frac{1}{2}\sigma(f)$.
Since $k\leq|V_{f}^{++}|$,
\begin{equation}
\sigma(f_{k})=\sigma(f)-k\geq\sigma(f)-|V_{f}^{++}|>\frac{1}{2}\sigma
(f).\label{eq_5.1b}%
\end{equation}

\item Assume therefore that $V_{f}^{1}=\varnothing$. In what follows, if
$k=1$, we take $f_{k-1}$ to be $f$ and ignore the reference to $v_{k-1}$. By
the construction of the broadcasts $f_{i}$, $i=1,...,k-1$, $\operatorname{PB}%
_{f_{i}}(v_{i})=B_{f_{i}}(v_{i})\neq\varnothing$. Indeed, we also have that
\begin{equation}
\operatorname{PB}_{f_{i}}(v_{j})=B_{f_{i}}(v_{j})=\operatorname{PB}_{f_{j}%
}(v_{j})\neq\varnothing\text{ for\ each\ }j\text{\ such\ that\ }1\leq j\leq
i.\label{eq_5.1}%
\end{equation}
Now, $\operatorname{PB}_{f_{k-1}}(v_{k})=\varnothing$ but $B_{f_{k-1}}%
(v_{k})\neq\varnothing$, hence there exists a vertex $u\in V_{f_{k-1}}^{+}$
such that $B_{f_{k-1}}(u)\cap B_{f_{k-1}}(v_{k})\neq\varnothing$. Since
$V_{f_{k-1}}^{+}=V_{f}^{+}$, $u\in V_{f}^{+}$, and since $V_{f}^{1}%
=\varnothing$, $u\in V_{f}^{++}$. Moreover, since $B_{f_{k-1}}(u)\cap
B_{f_{k-1}}(v_{k})\neq\varnothing$, $\operatorname{PB}_{f_{k-1}}(u)\neq
B_{f-1}(u)$. By (\ref{eq_5.1}), $u\notin\{v_{1},...,v_{k}\}$, that is, $u\in
V_{f}^{+}-\{v_{1},...,v_{k}\}$, and since $V_{f}^{+}=V_{f}^{++}$, we deduce
that $V_{f}^{++}-\{v_{1},...,v_{k}\}\neq\varnothing$. Hence $k<|V_{f}^{++}|$.
Similar to (\ref{eq_5.1b}),
\begin{equation}
\sigma(f_{k})=\sigma(f)-k>\sigma(f)-|V_{f}^{++}|\geq\frac{1}{2}\sigma
(f).\label{eq_5.1c}%
\end{equation}

\end{itemize}

Therefore, by (\ref{eq_5.1b}) and (\ref{eq_5.1c}),
\begin{equation}
\Gamma_{b}(G)\geq\sigma(f_{k})>\sigma(f)-\frac{1}{2}\sigma(f)=\frac{1}%
{2}\sigma(f)=\frac{1}{2}\alpha_{\operatorname{bn}}(G).\label{eq_5.1a}%
\end{equation}

Since we also have that $\alpha_{\operatorname{bnr}}(G)\leq\alpha
_{\operatorname{bn}}(G)$, it follows from (\ref{eq_5.1a}) that $\alpha
_{\operatorname{bnr}}(G)/\Gamma_{b}(G)\leq\alpha_{\operatorname{bn}}%
(G)/\Gamma_{b}(G)<2$.

{To show that the bound for }$\alpha_{\operatorname{bn}}(G)/\Gamma_{b}(G)$ is
asymptotically best possible, we consider the spider $S=\operatorname{Sp}%
(2^{k})$ for $k\geq3$. By Proposition \ref{Prop_Spider}$(i)$ and $(iii)$,
$\alpha_{\operatorname{bn}}(S)=2k$ and $\Gamma_{b}(S)=k+1$. Hence
\[
\lim_{k\rightarrow\infty}\alpha_{\operatorname{bn}}(S)/\Gamma_{b}%
(S)=\lim_{k\rightarrow\infty}2k/(k+1)=2.~\blacksquare
\]

\subsection{$\Gamma_{b}(G)/\alpha_{\operatorname{bn}(\operatorname{bnr})}(G)$}

Bouchemakh and Fergani \cite{BF} showed that if $G$ is a graph of order $n$
and minimum degree $\delta(G)$, then $\Gamma_{b}(G)\leq n-\delta(G)$. Since
$\alpha_{\operatorname{bnr}}(G)\geq\alpha(G)$, it follows that
\[
\frac{\Gamma_{b}(G)}{\alpha_{\operatorname{bn}}(G)}\leq\frac{\Gamma_{b}%
(G)}{\alpha_{\operatorname{bnr}}(G)}\leq\frac{n-\delta(G)}{\alpha(G)}.
\]
This leads to the following result.

\begin{proposition}
For any connected bipartite graph $G$, $\Gamma_{b}(G)/\alpha
_{\operatorname{bn}}(G)\leq\Gamma_{b}(G)/\alpha_{\operatorname{bnr}}(G)<2$.
\end{proposition}

\noindent\textbf{Proof.\hspace{0.1in}}Say $G$ has order $n$. Since $G$ is
bipartite, $\alpha(G)\geq\frac{n}{2}$. If $n=1$, the result is obvious. If
$n\geq2$, then%
\[
\frac{\Gamma_{b}(G)}{\alpha_{\operatorname{bn}}(G)}\leq\frac{\Gamma_{b}%
(G)}{\alpha_{\operatorname{bnr}}(G)}\leq\frac{n-\delta(G)}{\alpha(G)}\leq
\frac{2(n-1)}{n}<2.~\blacksquare
\]

\begin{theorem}
{For general graphs, the ratios }$\Gamma_{b}(G)/\alpha_{\operatorname{bnr}%
}(G)$ and $\Gamma_{b}(G)/\alpha_{\operatorname{bn}}(G)$ are unbounded.
\end{theorem}

\noindent\textbf{Proof.\hspace{0.1in}}Let $G_{n}\cong K_{n}%
\boksie
P_{3}$, where $X=\{x_{1},...,x_{n}\},\ Y=\{y_{1},...,y_{n}\}$ and
$Z=\{z_{1},...,z_{n}\}$ are the vertex sets of the copies of $K_{n}$, and
$Q_{i}=(x_{i},y_{i},z_{i})$ the copies of $P_{3}$. We begin by showing that
$\Gamma_{b}(G_{n})=2n$. If $n=1$, then $G_{1}=P_{3}$, hence $\Gamma_{b}%
(G_{1})=2$. Assume that $n\geq2$. Define the broadcast $f$ on $G_{n}$ by
$f(x)=2$ for $x\in X$ and $f(v)=0$ for $v\in Y\cup Z$. Then each $x_{i}$
broadcasts to all of $Q_{i}$ and $\operatorname{PB}_{f}(x_{i})=\{z_{i}\}$.
Hence $f$ is a minimal dominating broadcast. Consequently, $\Gamma_{b}%
(G_{n})\geq\sigma(f)=2n$.

Suppose there exists a minimal dominating broadcast $g$ on $G_{n}$ such that
$\sigma(g)>2n$. Then $n\geq2$. By the pigeonhole principle, there exists an
index $i$ such that $\sigma(Q_{i})>2$, and by symmetry we may assume $i=1$.
Since $\operatorname{diam}(G)=3$ and $e(x_{1})=e(z_{1})=3$ while $e(y_{1})=2$,
there are only three cases to consider.

\noindent\textbf{Case 1:\hspace{0.1in}}Without loss of generality,
$g(x_{1})=2$ and $g(z_{1})=1$. Since $N_{g}(x_{1})=X\cup Y\cup\{z_{1}\}$ and
$N_{g}(z_{1})=Z\cup\{y_{1}\}$, every vertex of $G_{n}$ hears $g$ from $x_{1}$
or $z_{1}$, hence $\sigma(g)=3$.\smallskip

\noindent\textbf{Case 2:\hspace{0.1in}}Without loss of generality,
$g(x_{1})=3$. Since $e(x_{1})=3$, $x_{1}$ dominates $G_{n}$ and $\sigma
(g)=3$.\smallskip

\noindent\textbf{Case 3:\hspace{0.1in}}Otherwise, $g(x_{1})=g(y_{1}%
)=g(z_{1})=1$. Then $\{x_{1},y_{1},z_{1}\}$ is a dominating set of $G_{n}$,
hence $\sigma(g)=3$.\smallskip

In each case we have a contradiction. We conclude that $\Gamma_{b}(G_{n})=2n$.

On the other hand, we show that $\alpha_{\operatorname{bn}}(G_{n})=3$ for each
$n\geq2$. Let $h$ be the characteristic function of a maximum independent set
$\{x_{i},y_{j},z_{k}\}$. Then $h$ is bn-independent, irredundant and
dominating. Note that $\{x_{j},y_{i}\}\subseteq B_{h}(x_{i})\cap B_{h}(y_{j})$
and $\{x_{k},z_{i}\}\subseteq B_{h}(x_{i})\cap B_{h}(z_{k})$. By Proposition
\ref{prop-max-bn}$(i)$, $h$ is a maximal bn-independent broadcast, thus
$\alpha_{\operatorname{bn}}(G_{n})\geq\alpha_{\operatorname{bnr}}(G)\geq
\sigma(h)=\alpha(G_{n})=3$. Consider any bn-independent broadcast $h^{\prime}$
on $G$. To maintain bn-independence, $h^{\prime}$ has at most one broadcasting
vertex in each of $X,\ Y$ and $Z$. Moreover, if a vertex $v$ broadcasts with
strength $2$, then either $v\in Y$ and dominates the entire graph, or, without
loss of generality, $v\in X$ and dominates $X$ and $Y$. In the latter case,
there is at most one other broadcasting vertex, say $z$, which belongs to $Z$,
and $h^{\prime}(z)=1$. Hence $\sigma(h^{\prime})\leq3$ and thus $\alpha
_{\operatorname{bn}}(G_{n})=\alpha_{\operatorname{bnr}}(G)=3$. It follows that
the {ratios }$\Gamma_{b}(G)/\alpha_{\operatorname{bnr}}(G)$ and $\Gamma
_{b}(G)/\alpha_{\operatorname{bn}}(G)$ are unbounded.~$\blacksquare$

\section{Open questions}

\label{Sec_Qs}As mentioned above, $\Gamma_{b}(G)\leq n-\delta(G)$ \cite{BF}
for all graphs $G$ of order $n$ and minimum degree $\delta(G)$. Also,
$\alpha(G)\leq n-\delta(G)$ and, when $G$ is connected, $\operatorname{diam}%
(G)\leq n-\delta(G)$.

\begin{question}
\label{Q_bn_upperbound}Is it true that $\alpha_{\operatorname{bn}}(G)\leq
n-\delta(G)$ for all graphs $G$?
\end{question}

In Proposition \ref{Prop_Gamma_minus_bn} we used $3\times n$ grids to show
that {$\Gamma_{b}-\alpha_{\operatorname{bn}}$ can be arbitrary for }%
$2$-connected graphs. Since $\alpha_{\operatorname{bn}}(T)\geq\alpha
_{\operatorname{bn}}(G)$ when $T$ is a spanning tree of $G$, this result does
not automatically extend to trees.

\begin{question}
\label{Q_bounded}{Is $\Gamma_{b}-\alpha_{\operatorname{bn}}$ bounded for
trees? }
\end{question}

We showed in Theorem \ref{gammabound} that $\alpha_{\operatorname{bnr}%
}(G)/\Gamma_{b}(G)<2$, but not that the bound is asymptotically best possible.
For the tree in Figure \ref{undombnr}, $\alpha_{\operatorname{bnr}}%
(T)/\Gamma_{b}(T)=14/13$, and for the general class of trees $H_{k}$
constructed using $T$, $\alpha_{\operatorname{bnr}}(H_{k})/\Gamma_{b}%
(H_{k})=42/41$. For the graph $G_{2}$ in Figure \ref{Fig_alpha_bnr},
$\alpha_{\operatorname{bnr}}(G_{2})/\Gamma_{b}(G_{2})=9/7$, and for the
general class of graphs $G_{k}$ constructed in Section \ref{Sec_bnr-Gam},
$\alpha_{\operatorname{bnr}}(G_{k})/\Gamma_{b}(G_{k})=(3k+3)/(2k+3)<3/2$.
Also, $\delta(G_{k})=1$.

\begin{question}
What is an asymptotically tight upper bound for $\alpha_{\operatorname{bnr}%
}(T)/\Gamma_{b}(T)$ for trees? Can the ratio $\alpha_{\operatorname{bnr}%
}(G)/\Gamma_{b}(G)<3/2$ for cyclic graphs be improved?
\end{question}

We showed in Theorem \ref{Thm_cyclic} that $\alpha_{\operatorname{bn}}%
-\Gamma_{b}$ and $\alpha${$_{\operatorname{bnr}}-\Gamma_{b}$ can be arbitrary
for cyclic graphs, but the graphs used in the proof have end-vertices.}

\begin{question}
\label{Q_2connected}Is it true that {$\alpha_{\operatorname{bnr}}(G)\leq
\Gamma_{b}(G)$ when }$G$ is{ }${2}${-connected? }Is it true that
{$\alpha_{\operatorname{bnr}}(G_{k})\leq\Gamma_{b}(G_{k})$ when }%
$\delta(G)\geq2$?
\end{question}

Dunbar et al.~\cite[Section 3.3]{Dunbar} also considered hearing independent
dominating broadcasts and denoted the maximum cost of a minimal independent
dominating broadcast of $G$, which they called the \emph{upper broadcast
independent domination number}, by $\Gamma_{\mathrm{ib}}(G)$. This parameter
was denoted $\alpha_{\mathrm{hd}}(G)$ in \cite{LindaD}. By Definition
\ref{bnd-i}, $\alpha_{\operatorname{bnd}}(G)\leq\alpha_{\mathrm{hd}}(G)$ for
all graphs $G$. The inequality can be strict: let $(u_{i},v_{i},w_{i}%
),\ i=1,2$, be two copies of $P_{3}$ and add the edges $v_{1}v_{2}$ and
$w_{1}w_{2}$ to form the graph $G$. The broadcast $f$ defined by $f(u_{i})=2$
for each $i$ is hearing independent and dominating. Since $w_{i}%
\in\operatorname{PB}_{f}(u_{i})$ for each $i$, it is minimal dominating. Hence
$\alpha_{\mathrm{hd}}(G)\geq4$. However, $\alpha_{\operatorname{bnd}%
}(G)=\alpha(G)=\operatorname{diam}(G)=3$. Bouchouika, Bouchemakh and Sopena
\cite{BBS} determined the values of the parameters defined in \cite{Dunbar}
for paths and cycles, but very little else is known about either of the
parameters $\alpha_{\operatorname{bnd}}$ and $\alpha_{\mathrm{hd}}$. We
encourage interested readers to investigate them.

\bigskip

\noindent\textbf{Acknowledgements\hspace{0.1in}}We acknowledge the support of
the Natural Sciences and Engineering Research Council of Canada (NSERC), PIN 253271.

\noindent Cette recherche a \'{e}t\'{e} financ\'{e}e par le Conseil de
recherches en sciences naturelles et en g\'{e}nie du Canada (CRSNG), PIN
253271.%
\begin{center}
\includegraphics[
natheight=0.773100in,
natwidth=1.599900in,
height=0.4151in,
width=0.8276in
]%
{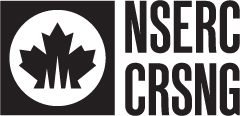}%
\end{center}


\begin{thebibliography}{99}                                                                                               %


\bibitem {Ahmadi}D.~Ahmadi, G.~H.~Fricke, C.~Schroeder, S.~T.~Hedetniemi and
R.~C.~Laskar, Broadcast irredundance in graphs. \emph{Congr.~Numer.}%
~\textbf{224} (2015), 17--31.

\bibitem {ABS}M.~Ahmane, I.~Bouchemakh and E.~Sopena, On the broadcast
independence number of caterpillars. \emph{Discrete Applied Math}.
\textbf{244} (2018), 20--356.

\bibitem {ABS2}M.~Ahmane, I.~Bouchemakh and E.~Sopena, On the broadcast
independence number of locally uniform 2-lobsters. arXiv:1902.02998v1, 2019.

\bibitem {BR}S.~Bessy and D.~Rautenbach, Relating broadcast independence and
independence. \emph{Discrete Math. }\textbf{342} (2019), 111589.
arXiv:1809.09288, 2018.

\bibitem {BR2}S.~Bessy and D.~Rautenbach, Girth, minimum degree, independence,
and broadcast independence. \emph{Commun. Comb. Optim.} \textbf{4} (2019), 131--139.

\bibitem {BF}I.~Bouchemakh and N.~Fergani, On the upper broadcast domination
number. \emph{Ars Combin.} \textbf{130} (2017), 151--161.

\bibitem {Bouch}I.~Bouchemakh and M.~Zemir, On the broadcast independence
number of grid graph. \emph{Graphs Combin.} \textbf{30} (2014), 83--100.

\bibitem {BBS}S.~Bouchouika, I.~Bouchemakh and E.~Sopena, Broadcasts on paths
and cycles. \emph{Discrete Appl. Math.} \textbf{283} (2020), 375--395.

\bibitem {CLZ}G.~Chartrand, L.~Lesniak and P.~Zhang, \emph{Graphs
\&\ Digraphs}, Chapman and Hall/CRC, Boca Raton, 2015.

\bibitem {Dunbar}J.~Dunbar, D.~Erwin, T.~Haynes, S.~M.~Hedetniemi and
S.~T.~Hedetniemi, Broadcasts in graphs. \emph{Discrete Applied Math.}
\textbf{154} (2006), 59-75.

\bibitem {Ethesis}D.~Erwin, \emph{Cost domination in graphs}. Doctoral
dissertation, Western Michigan University, 2001.

\bibitem {HMY}M.~A.~Henning, G.~MacGillivray and F.~Yang, Broadcast domination
in graphs. In T.~W.~Haynes, S.~T.~Hedetniemi and M.~A.~Henning (Eds.),
\emph{Structures of Domination in Graphs}, Springer, 2020, 15--46.

\bibitem {MM}E.~Marchessault and C.~M.~Mynhardt, Lower boundary independent
broadcasts in trees, submitted.

\bibitem {MN}C.~M.~Mynhardt, N.~Neilson, Boundary independent broadcasts in
graphs, \emph{J. Combin. Math. Combin. Comput. }\textbf{116} (2021), 79--100. https://doi.org/10.48550/arXiv.1906.10247

\bibitem {MN2}C.~M.~Mynhardt, N.~Neilson, A sharp upper bound for the boundary
independence broadcast number of a tree, submitted. https://doi.org/10.48550/arXiv.2104.02266

\bibitem {MR}C.~M.~Mynhardt, A.~Roux, Dominating and irredundant broadcasts in
graphs. \emph{Discrete Applied Math. }\textbf{220} (2017), 80-90.

\bibitem {LindaD}L.\ Neilson, \emph{Broadcast independence in graphs},
Doctoral dissertation, University of Victoria, 2019. http://hdl.handle.net/1828/11084
\end{thebibliography}
\end{document}